# A Complex-boundary Treatment Method for Finite Volume Schemes with Cut-Cartesian Cell Mesh

Zhaohui Qin

School of Engineering and Computer Science, Cedarville University, Cedarville, OH 45314

Abstract

In this paper, a second-order accurate method was developed for calculating fluid flows in complex geometries. This method uses cut-Cartesian cell mesh in finite volume framework. Calculus is employed to relate fluxes and gradients along curved surfaces to cell-averaged values. The resultant finite difference equations are sparse diagonal systems of equations. This method does not need repeated polynomial interpolation or reconstruction. Two-dimensional incompressible lid-driven semi-circular cavity flow at two Reynolds numbers was simulated with the current method and second-order accuracy was reached. The current method might be extended to third-order accuracy.

## 1. Introduction

Higher than second-order accurate numerical schemes are very desirable in many computational fluid dynamics (CFD) studies, for example in large eddy simulation (LES) of turbulent flows. In LES, the grid-scale flow motion is solved directly while the subgrid-scale (SGS) motion is modeled. Such a modeling practice usually introduces a SGS term proportional to the grid size squared into the momentum differential equations. The discretization error associated with a second-order scheme is of the same order as the SGS term, which interferes with the SGS term in an unpredictable way. Partly because of this reason, high-order accurate numerical schemes have become the norm in such studies. The spectral method (Orzag, 1980) and pseudo-spectral method (Ku, Hirsh, & Taylor, 1987), for instance, have been widely used in direct numerical simulation (DNS) and LES of turbulent flows, see, for example (Kim, Moin, & Moser, 1987). In the past few decades, high-order compact schemes (Lele, 1992) (Kobayashi, 1999) (Pereira, Kobayashi, & Pereira, 2001) based on Padé interpolation have gained increasing popularity in DNS and LES, for example, (Gamet, Ducros, Nicoud, & Poinsot, 1999) and (Xu L. , Cui, Xu, Zhang, & Chen, 2006). Compact schemes use shorter computational stencils than traditional finite difference and finite volume schemes yet yield the same order of accuracy (Lele, 1992) (Kobayashi, 1999).

It is unfortunately very difficult to apply high-order numerical schemes to flows with complex boundaries. For example, unstructured mesh has been extensively employed in simulations of flow with complex domains, but it is quite complicated to realize even second-order accuracy with such a mesh (Yu, Ozoe, & Tao, 2005). The spectral and pseudo-spectral methods are usually only effective for very simple flow geometries and boundary conditions. To apply compact schemes to flow problems involving irregular domains, special treatments must be developed to handle the complex boundary. (Piller & Stalio, 2008) employed a coordinate transformation method to recast compact schemes from the physical space to the computational space. (Lacor, Smirnov, & Baelmans, 2004) proposed a formulation which allows the use of compact schemes on arbitrary structured meshes. (Wang, Ren, & Li, 2016) devised a compact finite volume scheme on unstructured grids. (Xu L. , et al., 2006) implemented immersed boundary method (IBM) in conjunction with the compact scheme in their LES of flows with complex geometries. These methods usually rely on extensive polynomial fitting to reconstruct coefficients of the compact scheme at each irregular mesh cell, due to complexity of the cell geometry and connectivity. For example, in the work of (Lacor, Smirnov, & Baelmans, 2004), eight equations are solved at each irregular mesh cell to determine the coefficients of a third-order accurate scheme for two-dimensional inviscid flows. Such polynomial reconstructions must be carried out repeatedly as solution is being updated. It seems attractive to pursue high-order accurate schemes on a Cartesian mesh, which may enable us to avoid the laborious polynomial reconstructions due to the exceeding simplicity of such a mesh. In fact, one of the earliest successes in CFD was accomplished by adopting a cut-Cartesian cell mesh and finite difference method to solve the flow around a cylinder (Thom, 1933). It is, however, challenging to incorporate the cut-Cartesian cell mesh in the finite volume context. (Ye, Mittal, Udaykumar, & Shyy, 1999) developed a polynomial interpolation method to obtain a second-order accurate approximation of the fluxes and gradients on the faces of Cartesian cells cut by immersed boundaries. Six unknown coefficients of a polynomial must be solved at each cut-cell face from available neighboring cell-center values, and on top of which sixteen different cut-cell configurations have to be considered for two-dimensional flows. Such polynomial interpolations take place at every time step. (Mariani & Prata, 2008) employed this method to explore fluid flows in two-dimensional lid-driven irregular-bottom cavities.

Another issue that permeates most finite volume-based complex-boundary treatment methods is the approximation of curved boundaries with piecewise line segments, which is only first-order accurate in terms of geometry reproduction precision (Pilliod & Puckett, 2004).

The objective of the current paper is to present a complex-boundary treatment method for finite volume schemes using cut-Cartesian cell mesh. As the first step toward a high-order finite volume method, we pursue a second-order accurate scheme. Two-dimensional incompressible lid-driven semi-circular cavity flow is selected as an example to illustrate the effectiveness of the current method.

Consider a two-dimensional semi-circular cavity of diameter $D$ as shown in Figure 1. The fluid in the cavity is of constant density and constant viscosity. The fluid flow is driven by the cavity lid, which is infinitely long and translates to the right with a constant speed $U$. This flow has been studied by (Glowinski, Guidoboni, & Pan, 2006) in a benchmark work with an operator-splitting/finite elements-based method. Their results and data obtained using a commercial CFD software with a very fine mesh will be utilized to validate the current method.

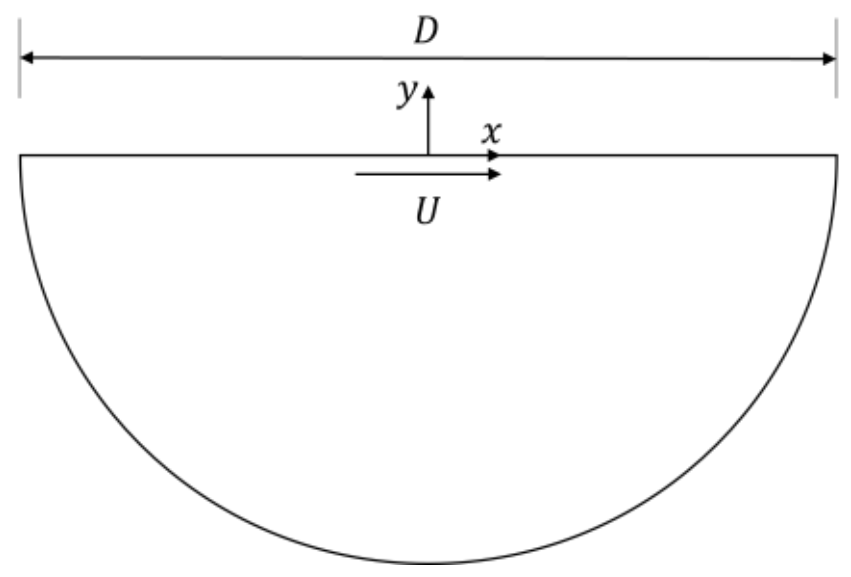

Figure 1

As will be demonstrated, the current method does not need repeated polynomial reconstructions or interpolations. It also does not approximate irregular boundaries with piecewise linear segments. Moreover, this method permits enormous flexibility in technical details when it is applied. The paper is organized as follows. The numerical method is expounded in Section 2. In Section 3, numerical results of the lid-driven semi-circular cavity flow are presented and compared with benchmark data. The paper is concluded with a summary of the current method and a discussion of the possible route to higher-order schemes.

# 2. Numerical Method

## 2.1.Governing differential equations

The flow field in the cavity is governed by the Navier-Stokes and continuity equations

$$\frac{\partial u}{\partial t}+\frac{\partial (uu)}{\partial x}+\frac{\partial (vu)}{\partial y}=-\frac{\partial p}{\partial x}+\nu\left(\frac{\partial^2 u}{\partial x^2}+\frac{\partial^2 u}{\partial y^2}\right) \tag{1}$$

$$\frac{\partial v}{\partial t}+\frac{\partial (uv)}{\partial x}+\frac{\partial (vv)}{\partial y}=-\frac{\partial p}{\partial y}+\nu\left(\frac{\partial^2 v}{\partial x^2}+\frac{\partial^2 v}{\partial y^2}\right) \tag{2}$$

$$\frac{\partial u}{\partial x}+\frac{\partial v}{\partial y}=0 \tag{3}$$

where $t$ is time; $u$ and $v$ are flow velocities along the $x$- and $y$-directions; $p$ is the kinematic pressure and $\nu$ is the kinematic viscosity.

No-slip conditions apply at the boundaries, i.e.,

$$u=v=0 \tag{4}$$

along the semi-circular wall and

$$u=U, v=0 \tag{5}$$

at the cavity lid.

## 2.2.Mesh

A Cartesian mesh is employed to discretize the governing differential equations. Cut cells are formed as Cartesian cells are cut by boundaries as depicted in Figure 2. This mesh, however, needs a few modifications before it may satisfy three requirements demanded by the current method.

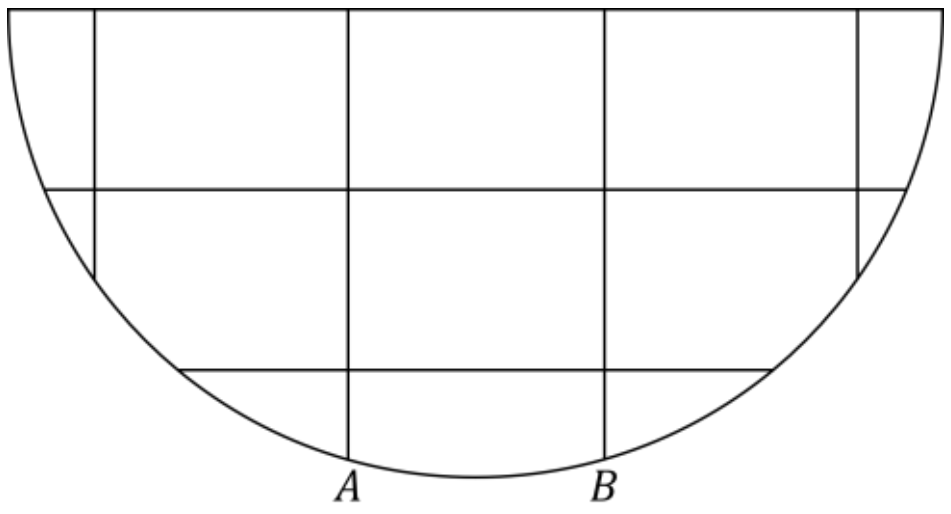

Figure 2

First, each cell face must be smooth, i.e., without discontinuities in the face curvature. Secondly, each coordinate of a cell face must be monotonic. Thirdly, two grid lines, one horizontal and one vertical, must be ruled at each end of a cell face. The rationale behind these requirements will be clear later on.

The mesh shown in Figure 2 does not meet the second requirement as the $y$-coordinate of face AB is not monotonic. It also fails to fulfill the third requirement, for example there is only one vertical grid line at each end of face AB.

The first two requirements can be easily satisfied by properly splitting the geometry boundaries when we mesh the geometry. The third requirement, however, is only feasible for computational domains with at least one straight boundary that spans the whole domain along that direction (e.g., the cavity lid) or symmetrical domains, which the cavity also happens to be. Otherwise, enforcement of this requirement may result in infinite grid lines. Fortunately, most practical computational domains either belong to or can be easily adapted to such shapes without compromising the flow physics.

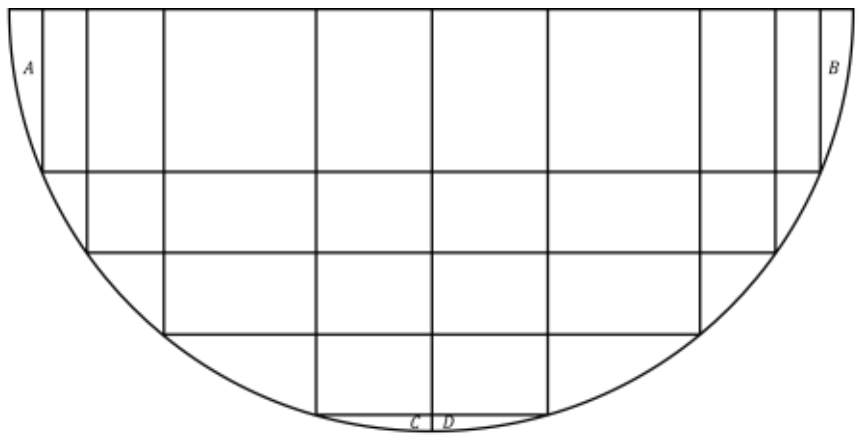

Figure 3

A mesh complying with all three requirements is shown in Figure 3. As can be observed, such a mesh contains only four-faced rectangular cells and three-faced cut cells. The four cut cells marked with letters call for special treatments. Cells A and B are termed solitary cells as they do not have neighboring cells along one direction ($y$-direction here). Cells C and D are named twin cells, which are side-by-side cut cells. Solitary cells and twin cells are the only outliers in any two-dimensional mesh to which the current method may apply.

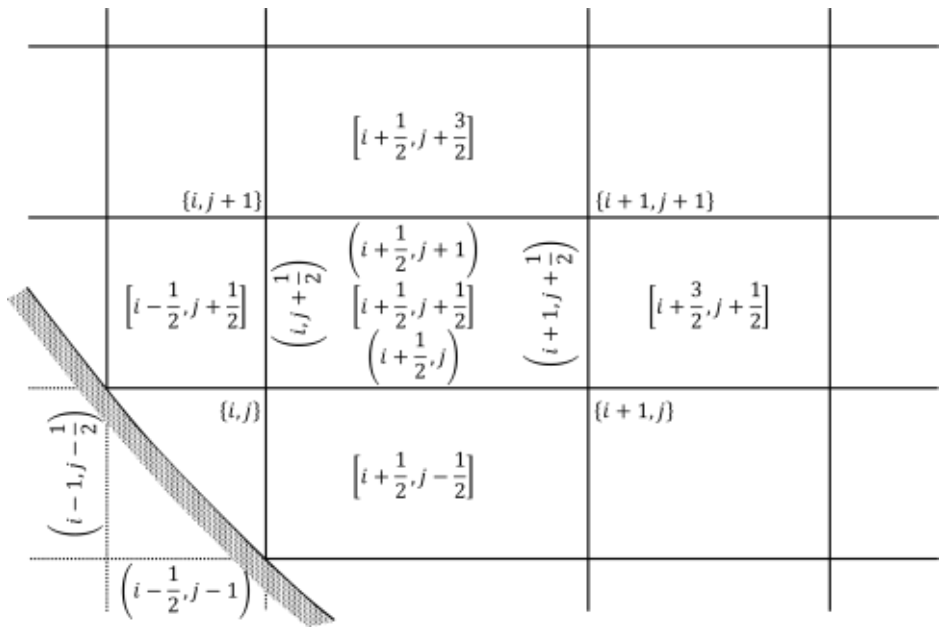

*Figure 4*

In the present work we number the "original", i.e., the pre-cut Cartesian mesh first. A cut cell then inherits the index of the original pre-cut cell, and a curved boundary face inherits both indices of the two faces it cuts off. For example, the complete cut cell at the lower left corner of Figure 4 is the $\left[i-\frac{1}{2},j-\frac{1}{2}\right]$ cell, and the curved boundary face of this cell is both the $\left(i-1,j-\frac{1}{2}\right)$ face and the $\left(i-\frac{1}{2},j-1\right)$ face. In this figure indices in square brackets are of cells; those in parentheses are of faces and those in curly brackets are of nodes. In subsequent sections a single index is sometimes used to denote a mesh entity if confusion is not caused.

## 2.3. Discretization

The governing differential equations are discretized by integrating them over each mesh cell. Such integrations can be done in two hypothetical steps. In the first step, the equation is integrated over a thin slice of the cell. For example, for cell $\left[i-\frac{1}{2},j-\frac{1}{2}\right]$, one may cut a slice of it with a horizontal line $y=Y$, as shown in Figure 5. Integrating equation (1) over this slice results in

$$\begin{aligned}\frac{\partial}{\partial t}u^X dy+[(uu)_e-(uu)_w]dy+\left[\frac{\partial(vu)}{\partial y}\right]^X dy \\ =-(p_e-p_w)dy+\nu\left[\left(\frac{\partial u}{\partial x}\right)_e-\left(\frac{\partial u}{\partial x}\right)_w+\left(\frac{\partial u}{\partial y}\right)^X\right]dy\end{aligned} \tag{6}$$

in which subscripts $w$ and $e$ denote values at the west (left) and east (right) ends of the slice. Here we define slice integral of a function $\phi(x,y)$ as

$$\phi^X(Y)=\int_{x_w(Y)}^{x_e(Y)}\phi(x,Y)dx \tag{7}$$

where $x_w(Y)$ and $x_e(Y)$ are the $x$-cooredinates of the two ends of the cell slice cut out by $y=Y$. In the second step, Equation (6) is integrated over the whole $y$-range of the cell. Alternatively, one may integrate the governing differential equation over a vertical cell slice first, then over the whole cell $x$-range. If this order of integrations is preferred, one may define the slice integral along a vertical cell slice cut out by line $x=X$ as

$$\phi^{Y}(X)=\int_{y_s(X)}^{y_n(X)}\phi(X,y)dy \tag{8}$$

where $y_s(X)$ and $y_n(X)$ are the $y$-cooredinates of the south (bottom) and north (top) ends of the slice. Following these two steps, which integration order we may choose does not matter, we have

$$\begin{aligned}\frac{\partial}{\partial t}u^{xy}_{\left[i-\frac{1}{2},j-\frac{1}{2}\right]}&+\left[(uu)^{y}_{\left(i,j-\frac{1}{2}\right)}-(uu)^{y}_{\left(i-1,j-\frac{1}{2}\right)}\right]+\left[(vu)^{x}_{\left(i-\frac{1}{2},j\right)}-(vu)^{x}_{\left(i-\frac{1}{2},j-1\right)}\right]\\ &=-\left[p^{y}_{\left(i,j-\frac{1}{2}\right)}-p^{y}_{\left(i-1,j-\frac{1}{2}\right)}\right]\\ &+\nu\left[\left(\frac{\partial u}{\partial x}\right)^{y}_{\left(i,j-\frac{1}{2}\right)}-\left(\frac{\partial u}{\partial x}\right)^{y}_{\left(i-1,j-\frac{1}{2}\right)}+\left(\frac{\partial u}{\partial y}\right)^{x}_{\left(i-\frac{1}{2},j\right)}\right.\\ &\left.-\left(\frac{\partial u}{\partial y}\right)^{x}_{\left(i-\frac{1}{2},j-1\right)}\right]\end{aligned} \tag{9}$$

where the cell integral is

$$\phi^{xy}=\int_{x_w}^{x_e}\phi^{Y}(X)dX=\int_{y_s}^{y_n}\phi^{X}(Y)dY \tag{10}$$

in which $x_w$, $x_e$, $y_s$ and $y_n$ are the lower and upper coordinate limits of the cell in the $x$- and $y$-directions. We define the sliding integrals of a function $\phi$ along a face, say face $f$ as

$$\phi^{x}_{f}=\int_{x_w}^{x_e}\phi_f dx \tag{11}$$

$$\phi^{y}_{f}=\int_{y_s}^{y_n}\phi_f dy \tag{12}$$

where $\phi_f$ is the $\phi$ value along face $f$. Notice that both sliding integrals along an inclined face may be nonzero. With the same procedure, Equation (2) can be integrated and becomes

$$\begin{aligned}\frac{\partial}{\partial t}v^{xy}_{\left[i-\frac{1}{2},j-\frac{1}{2}\right]}&+\left[(uv)^{y}_{\left(i,j-\frac{1}{2}\right)}-(uv)^{y}_{\left(i-1,j-\frac{1}{2}\right)}\right]+\left[(vv)^{x}_{\left(i-\frac{1}{2},j\right)}-(vv)^{x}_{\left(i-\frac{1}{2},j-1\right)}\right]\\ &=-\left[p^{x}_{\left(i-\frac{1}{2},j\right)}-p^{x}_{\left(i-\frac{1}{2},j-1\right)}\right]\\ &+\nu\left[\left(\frac{\partial v}{\partial x}\right)^{y}_{\left(i,j-\frac{1}{2}\right)}-\left(\frac{\partial v}{\partial x}\right)^{y}_{\left(i-1,j-\frac{1}{2}\right)}+\left(\frac{\partial v}{\partial y}\right)^{x}_{\left(i-\frac{1}{2},j\right)}\right.\\ &\left.-\left(\frac{\partial v}{\partial y}\right)^{x}_{\left(i-\frac{1}{2},j-1\right)}\right]\end{aligned} \tag{13}$$

Notice that equations (9) and (13) are exact as no assumptions have been made yet, and they are valid for both interior rectangular cells and boundary cut cells.

It is convenient to define the slice averages, cell average, and sliding averages as follows:

$$\bar{\phi}^X(Y) = \frac{\phi^X(Y)}{x_e(Y) - x_w(Y)} \tag{14}$$

$$\bar{\phi}^Y(X) = \frac{\phi^Y(X)}{y_n(X) - y_s(X)} \tag{15}$$

$$\bar{\phi}^{xy} = \frac{\phi^{xy}}{\int_{x_w}^{x_e}[y_n(X) - y_s(X)]dX} = \frac{\phi^{xy}}{\int_{y_s}^{y_n}[x_e(Y) - x_w(Y)]dY} = \frac{\phi^{xy}}{\Delta V} \tag{16}$$

$$\bar{\phi}_f^x = \frac{\phi_f^x}{x_e - x_w} = \frac{\phi_f^x}{\Delta x} \tag{17}$$

$$\bar{\phi}_f^y = \frac{\phi_f^y}{y_n - y_s} = \frac{\phi_f^y}{\Delta y} \tag{18}$$

in which $\Delta x$, $\Delta y$ and $\Delta V$ are the overall $x$-, $y$-sizes and volume of the cell. One should notice that as a slice shrinks to a node point, the slice average is simply the nodal value. For example, for cell $\left[i - \frac{1}{2}, j - \frac{1}{2}\right]$ depicted in Figure 5, the slice average along $y = y_s$ is $\bar{\phi}^X(y_s) = \phi_{\{i,j-1\}}$. We then define the mean $x$- and $y$-sizes of a cell as

$$\overline{\Delta x} = \frac{\int_{y_s}^{y_n}[x_e(Y) - x_w(Y)]dY}{y_n - y_s} = \frac{\Delta V}{\Delta y} \tag{19}$$

$$\overline{\Delta y} = \frac{\int_{x_w}^{x_e}[y_n(X) - y_s(X)]dX}{x_e - x_w} = \frac{\Delta V}{\Delta x} \tag{20}$$

Notice that $\Delta x$, $\Delta y$, $\Delta V$, $\overline{\Delta x}$, and $\overline{\Delta y}$ are all constants, which only need to be calculated once. The integrated momentum equations (9) and (13) can then be rewritten as

$$\frac{\partial}{\partial t}\bar{u}^{xy}_{\left[i-\frac{1}{2},j-\frac{1}{2}\right]}+\frac{\overline{(uu)}^{y}_{\left(i,j-\frac{1}{2}\right)}-\overline{(uu)}^{y}_{\left(i-1,j-\frac{1}{2}\right)}}{\overline{\Delta x}}+\frac{\overline{(vu)}^{x}_{\left(i-\frac{1}{2},j\right)}-\overline{(vu)}^{x}_{\left(i-\frac{1}{2},j-1\right)}}{\overline{\Delta y}}$$
$$=-\frac{\bar{p}^{y}_{\left(i,j-\frac{1}{2}\right)}-\bar{p}^{y}_{\left(i-1,j-\frac{1}{2}\right)}}{\overline{\Delta x}} \tag{21}$$
$$+\nu\left[\frac{\overline{\left(\frac{\partial u}{\partial x}\right)}^{y}_{\left(i,j-\frac{1}{2}\right)}-\overline{\left(\frac{\partial u}{\partial x}\right)}^{y}_{\left(i-1,j-\frac{1}{2}\right)}}{\overline{\Delta x}}+\frac{\overline{\left(\frac{\partial u}{\partial y}\right)}^{x}_{\left(i-\frac{1}{2},j\right)}-\overline{\left(\frac{\partial u}{\partial y}\right)}^{x}_{\left(i-\frac{1}{2},j-1\right)}}{\overline{\Delta y}}\right]$$

and

$$\frac{\partial}{\partial t}\bar{v}^{xy}_{\left[i-\frac{1}{2},j-\frac{1}{2}\right]}+\frac{\overline{(uv)}^{y}_{\left(i,j-\frac{1}{2}\right)}-\overline{(uv)}^{y}_{\left(i-1,j-\frac{1}{2}\right)}}{\overline{\Delta x}}+\frac{\overline{(vv)}^{x}_{\left(i-\frac{1}{2},j\right)}-\overline{(vv)}^{x}_{\left(i-\frac{1}{2},j-1\right)}}{\overline{\Delta y}}$$
$$=-\frac{\bar{p}^{x}_{\left(i-\frac{1}{2},j\right)}-\bar{p}^{x}_{\left(i-\frac{1}{2},j-1\right)}}{\overline{\Delta y}} \tag{22}$$
$$+\nu\left[\frac{\overline{\left(\frac{\partial v}{\partial x}\right)}^{y}_{\left(i,j-\frac{1}{2}\right)}-\overline{\left(\frac{\partial v}{\partial x}\right)}^{y}_{\left(i-1,j-\frac{1}{2}\right)}}{\overline{\Delta x}}+\frac{\overline{\left(\frac{\partial v}{\partial y}\right)}^{x}_{\left(i-\frac{1}{2},j\right)}-\overline{\left(\frac{\partial v}{\partial y}\right)}^{x}_{\left(i-\frac{1}{2},j-1\right)}}{\overline{\Delta y}}\right]$$

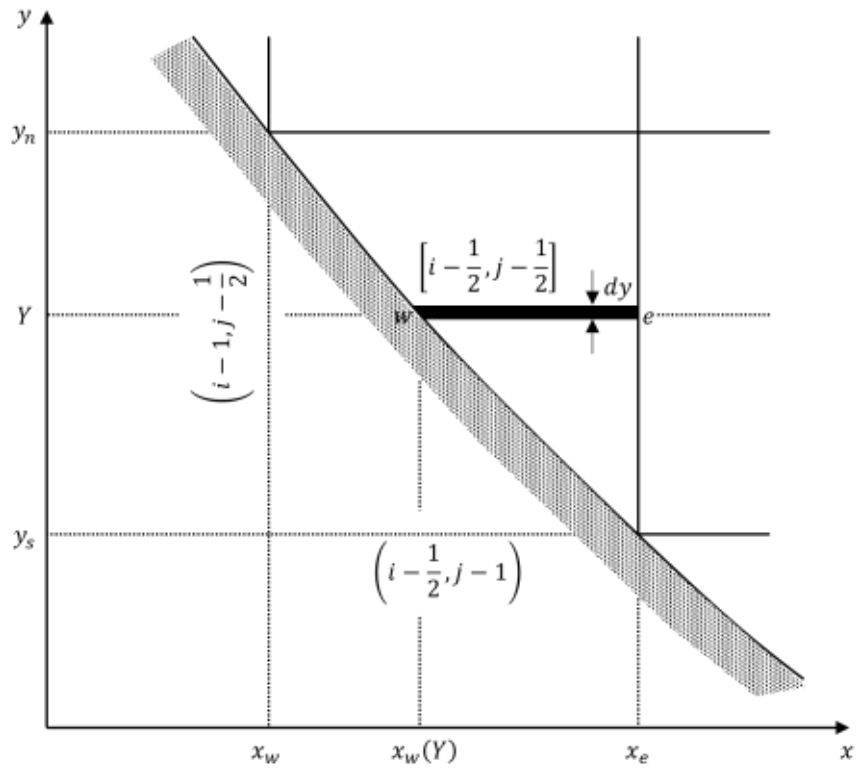

*Figure 5*

The objective is to obtain cell averages of the velocity field to the second-order spatial accuracy, which requires to evaluate sliding averages in equations (21) and (22) with second-order truncation errors (Kobayashi, 1999). Moreover, cell averages are the variables to be solved, therefore sliding averages must be related to corresponding cell averages if they are unknown from boundary conditions. To this end, we will first introduce the PKP theorem, which is of fundamental importance in the development of the current method.

## 2.4. PKP theorem

If two functions $\phi$ and $\psi$ can be regarded as functions of only one coordinate, say $y$, by comparing the Taylor expansion of the integral of the product of these two functions with that of the product of integrals of these two functions, we have (Pereira, Kobayashi, & Pereira, 2001)

$$\frac{\int_{y_s}^{y_n}(\phi\psi)dy}{\Delta y}=\frac{\int_{y_s}^{y_n}\phi dy}{\Delta y}\frac{\int_{y_s}^{y_n}\psi dy}{\Delta y}+\frac{\Delta y}{12}\left(\frac{d\phi}{dy}\right)_0\left(\frac{d\psi}{dy}\right)_0+\mathcal{O}(\Delta y^4) \tag{23}$$

in which $\Delta y=y_n-y_s$ and the subscript 0 denotes the middle point of this $y$-range. Equation (24) will be termed PKP theorem (after its authors) hitherto. Clearly this theorem is applicable to sliding averages since variables along a face can be treated as functions of only one coordinate. Therefore,

$$\overline{(\phi\psi)}^{y}_{f}=\bar{\phi}^{y}_{f}\bar{\psi}^{y}_{f}+\frac{\Delta y_f^2}{12}\left(\frac{d\phi_f}{dy}\right)_0\left(\frac{d\psi_f}{dy}\right)_0+\mathcal{O}\left(\Delta y_f^4\right) \tag{24}$$

Notice that the two derivatives on the right side of this equation are total derivatives rather than partial derivatives. For example, the total $y$-derivative of function $\phi$ at point $b$ along the face sketched in Figure 6 is

$$\left(\frac{d\phi_f}{dy}\right)_b = \lim_{\delta y \to 0} \frac{\phi_a - \phi_b}{\delta y} \tag{25}$$

On the other hand, the partial $y$-derivative of $\phi$ at the same point is

$$\left(\frac{\partial\phi}{\partial y}\right)_b = \lim_{\delta y \to 0} \frac{\phi_c - \phi_b}{\delta y} \tag{26}$$

Clearly, they coincide along a vertical face, but they may differ greatly for an almost-horizontal inclined face. These two derivatives are related by

$$\left(\frac{d\phi_f}{dy}\right)_b = \left(\frac{\partial\phi}{\partial y}\right)_b + \left(\frac{\partial\phi}{\partial x}\right)_b \left(\frac{dx_f}{dy_f}\right)_b \tag{27}$$

in which $\left(dx_f/dy_f\right)_b$ is the inversed slope of face $f$ at point $b$. Evidently

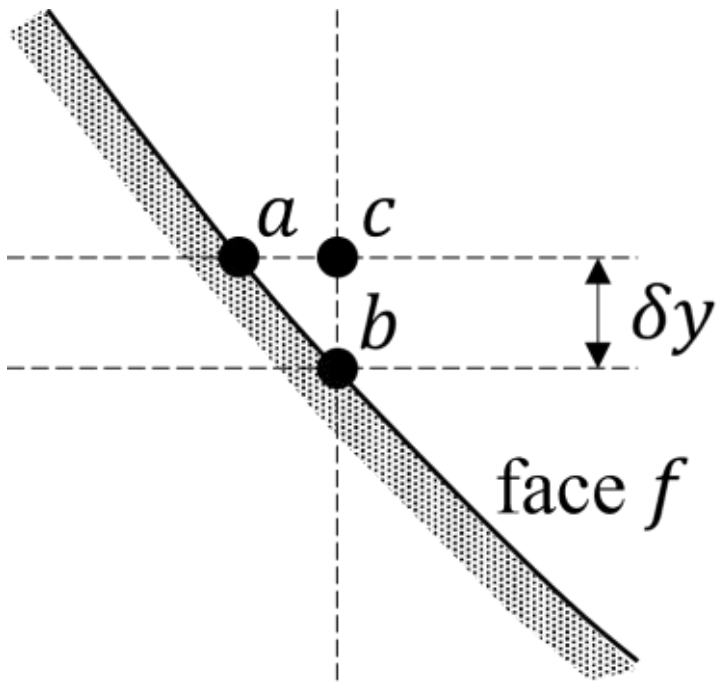

Figure 6

$$\left(\frac{d\phi_f}{dy}\right)_b = \left(\frac{d\phi_f}{dx}\right)_b \left(\frac{dx_f}{dy_f}\right)_b \tag{28}$$

and

$$\left(\frac{d\phi_f}{dx}\right)_b = \left(\frac{d\phi_f}{dy}\right)_b \left(\frac{dy_f}{dx_f}\right)_b \tag{29}$$

If $\overline{(\phi\psi)}_f^y$ is to be just second-order accurate, for example the nonlinear terms in equations (21) and (22), only the first term on the right side of Equation (24) is needed. However, to reach third-order accuracy in $\overline{(\phi\psi)}_f^y$ as is necessary in later development of the current method, the second term on the right side of Equation (24) must be retained and the two total derivatives in this term must be approximated with at least first-order accuracy. This can be done in many ways. For example, we may use sliding averages along the south and north neighboring faces of a vertical or inclined face to approximate the total derivative at the $y$-center of such a face:

$$\left(\frac{d\phi_f}{dy}\right)_0 = a_s\bar{\phi}_s^y + a_n\bar{\phi}_n^y + \mathcal{O}(h) \tag{30}$$

where

$$\begin{cases} a_s = -\dfrac{2}{2\Delta y_f + \Delta y_s + \Delta y_n} \\ a_n = \dfrac{2}{2\Delta y_f + \Delta y_s + \Delta y_n} \end{cases} \tag{31}$$

and subscripts $s$ and $n$ refer to the immediate south and north neighboring faces. The symbol $h$ represents a constant comparable to the local mesh size. This equation can be derived by using Taylor expansion. Notice that the south/north neighboring face of an inclined boundary face is the boundary face right below/above it. For example, the south neighboring face of boundary face $\left(i-1, j-\frac{1}{2}\right)$ in Figure 4 is boundary face $\left(i, j-\frac{3}{2}\right)$.

If face $f$ has only north neighboring faces, we have

$$\left(\frac{d\phi_f}{dy}\right)_0 = a_f\bar{\phi}_f^y + a_n\bar{\phi}_n^y + \mathcal{O}(h) \tag{32}$$

where

$$\begin{cases} a_f = -\dfrac{2}{\Delta y_f + \Delta y_n} \\ a_n = \dfrac{2}{\Delta y_f + \Delta y_n} \end{cases} \tag{33}$$

On the other hand, if face $f$ has only south neighboring faces, we have

$$\left(\frac{d\phi_f}{dy}\right)_0 = a_s\bar{\phi}_s^y + a_f\bar{\phi}_f^y + \mathcal{O}(h) \tag{34}$$

where

$$\begin{cases} a_s = -\dfrac{2}{\Delta y_s + \Delta y_f} \\ a_f = \dfrac{2}{\Delta y_s + \Delta y_f} \end{cases} \tag{35}$$

Alternatively, we can use the familiar central difference to approximate such derivatives:

$$\left(\frac{d\phi_f}{dy}\right)_0 = \frac{\phi_f(y_n) - \phi_f(y_s)}{\Delta y_f} + \mathcal{O}(h^2) \tag{36}$$

where $y_s$ and $y_n$ are the lower and upper $y$-limits of the face. This approximation has at least first-order accuracy if the $\phi$ values at the two ends of the face are second-order accurate. Such nodal values can be obtained from sliding averages of neighboring faces, as follows:

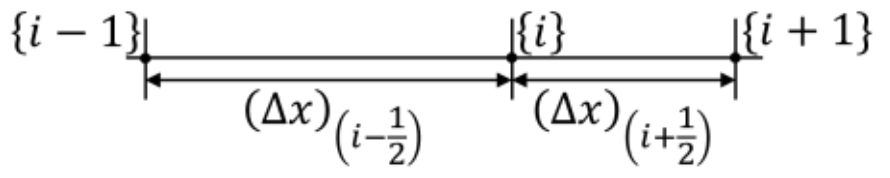

*Figure 7*

For three consecutive nodes, say nodes $\{i-1\}$, $\{i\}$, and $\{i+1\}$ shown in Figure 7, we have

$$b_i\phi_{\{i\}} = d_i\bar{\phi}^x_{\left(i-\frac{1}{2}\right)} + e_i\bar{\phi}^x_{\left(i+\frac{1}{2}\right)} + \mathcal{O}(h^3) \tag{37}$$

The coefficients are

$$\begin{cases} b_i = (\Delta x)_{\left(i-\frac{1}{2}\right)} + (\Delta x)_{\left(i+\frac{1}{2}\right)} \\ d_i = (\Delta x)_{\left(i+\frac{1}{2}\right)} \\ e_i = (\Delta x)_{\left(i-\frac{1}{2}\right)} \end{cases} \tag{38}$$

This formula can be derived by equating the corresponding terms of the Taylor's series of both sides of Equation (37) up to the desired order of accuracy.

If node $\{i\}$ is located on a west (left) boundary and if $\phi$ is unknown along this boundary, it can be evaluated by the following formula

$$\phi_{\{i\}} + \phi_{\{i+1\}} = 2\bar{\phi}^x_{\left(i+\frac{1}{2}\right)} + \mathcal{O}(h^2) \tag{39}$$

Here a compact scheme is adopted to minimize the stencil size.

If node $\{i\}$ is on an east (right) boundary, we have

$$\phi_{\{i-1\}} + \phi_{\{i\}} = 2\bar{\phi}^x_{\left(i-\frac{1}{2}\right)} + \mathcal{O}(h^2) \tag{40}$$

One should notice that, as mentioned in the introduction of this paper, great flexibility is allowed in application of the current method. For example, an interior nodal value may be calculated from the sliding averages along neighboring vertical instead of horizontal faces. Also, one may directly relate a nodal value to the surrounding cell averages instead of sliding averages. Another example is when one evaluates a $y$-derivative at the center of a boundary face, she may use Equation (32) or (33); alternatively, she can compute the $x$-derivative first then use Equation (28) to obtain the $y$-derivative. Moreover, one has the freedom to choose any numerical scheme at almost every step of the procedure, as long as the desired accuracy is maintained. Now we may turn our attention to relating sliding averages to cell averages.

### 2.5. Evaluation of sliding averages

Sliding averages, if not obtainable from boundary conditions, can be related to cell averages as follows.

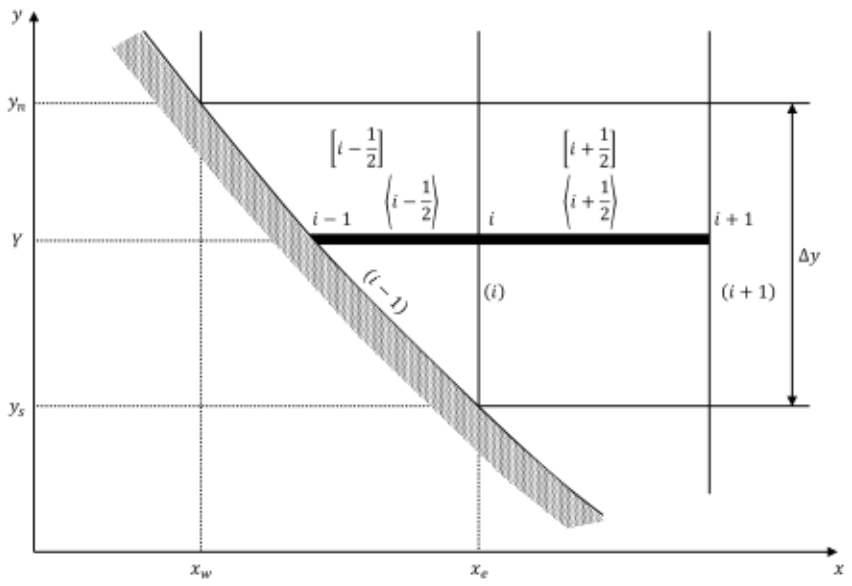

*Figure 8*

We may use a horizontal line $y = Y$ to cut the two cells straddling the face of interest, say face $(i)$ as shown in Figure 8. We have a formula similar to Equation (37):

$$b_i \phi_i = d_i \bar{\phi}^X_{\langle i-\frac{1}{2}\rangle} + e_i \bar{\phi}^X_{\langle i+\frac{1}{2}\rangle} + \mathcal{O}(h^3) \tag{41}$$

The coefficients are

$$\begin{cases} b_i = (\Delta x)_{\langle i-\frac{1}{2}\rangle} + (\Delta x)_{\langle i+\frac{1}{2}\rangle} \\ d_i = (\Delta x)_{\langle i+\frac{1}{2}\rangle} \\ e_i = (\Delta x)_{\langle i-\frac{1}{2}\rangle} \end{cases} \tag{42}$$

where subscripts in angle brackets denote slices in corresponding cells. Here $\phi_i$ is the $\phi$ value at the intersection of the horizontal line and face $(i)$. Integrating equation (41) over the $y$-range of face $(i)$, then dividing it by $\Delta y$, we have

$$\frac{\int_{y_s}^{y_n} (b_i \phi_i) dy}{\Delta y} = \frac{\int_{y_s}^{y_n} \left( d_i \bar{\phi}^X_{\langle i-\frac{1}{2}\rangle} \right) dy}{\Delta y} + \frac{\int_{y_s}^{y_n} \left( e_i \bar{\phi}^X_{\langle i+\frac{1}{2}\rangle} \right) dy}{\Delta y} + \mathcal{O}(h^3) \tag{43}$$

If the two cells on both sides of face $(i)$ are both rectangular cells, coefficients $b_i$, $d_i$, $e_i$, and $(\Delta x)_{<i\pm\frac{1}{2}>}$ are constants and equation (43) reduces to

$$b_i \bar{\phi}^y_{(i)} = d_i \bar{\phi}^{xy}_{\left[i-\frac{1}{2}\right]} + e_i \bar{\phi}^{xy}_{\left[i+\frac{1}{2}\right]} + \mathcal{O}(h^3) \tag{44}$$

If at least one of these cells is an irregular boundary cell, coefficients $b_i$, $d_i$, $e_i$ and $(\Delta x)_{<i\pm\frac{1}{2}>}$ might no longer be constants but known functions of the $y$-coordinate. Terms in Equation (43) then can be evaluated with the PKP theorem. The left side of Equation (43) is

$$\frac{\int_{y_s}^{y_n} (b_i \phi_i) dy}{\Delta y} = \frac{\int_{y_s}^{y_n} b_i dy}{\Delta y} \bar{\phi}^y_{(i)} + \frac{\Delta y^2}{12} \left( \frac{db_i}{dy} \right)_0 \left[ \frac{d\phi_{(i)}}{dy} \right]_0 + \mathcal{O}(h^4) \tag{45}$$

Since third-order accuracy is needed for this term, both terms on the right side of Equation (45) must remain and the two derivatives on this side must be at least first-order accurate. Evidently

$$\frac{\int_{y_s}^{y_n} b_i dy}{\Delta y} = \overline{(\Delta x)}_{\left[i-\frac{1}{2}\right]} + \overline{(\Delta x)}_{\left[i+\frac{1}{2}\right]} \tag{46}$$

and

$$\left(\frac{db_i}{dy}\right)_0 = \left[\frac{d(\Delta x)_{\langle i-\frac{1}{2}\rangle}}{dy}\right]_0 + \left[\frac{d(\Delta x)_{\langle i+\frac{1}{2}\rangle}}{dy}\right]_0 \tag{47}$$

can be calculated exactly and only need to be done once. The derivative $[d\phi_{(i)}/dy]_0$ can be assessed with the methods discussed in the previous section.

Similarly, the right side of Equation (43) is

$$\begin{aligned}\frac{\int_{y_s}^{y_n}\left(d_i\bar{\phi}^X_{\langle i-\frac{1}{2}\rangle}\right)dy}{\Delta y} + \frac{\int_{y_s}^{y_n}\left(e_i\bar{\phi}^X_{\langle i+\frac{1}{2}\rangle}\right)dy}{\Delta y} &= \overline{(\Delta x)}_{\left[i+\frac{1}{2}\right]}\overline{\left(\bar{\phi}^X\right)}^y_{\left[i-\frac{1}{2}\right]} + \overline{(\Delta x)}_{\left[i-\frac{1}{2}\right]}\overline{\left(\bar{\phi}^X\right)}^y_{\left[i+\frac{1}{2}\right]} \\ &+ \frac{\Delta y^2}{12}\left\{\left[\frac{d(\Delta x)_{\langle i+\frac{1}{2}\rangle}}{dy}\right]_0\left[\frac{d\bar{\phi}^X_{\langle i-\frac{1}{2}\rangle}}{dy}\right]_0 \right. \\ &\left. + \left[\frac{d(\Delta x)_{\langle i-\frac{1}{2}\rangle}}{dy}\right]_0\left[\frac{d\bar{\phi}^X_{\langle i+\frac{1}{2}\rangle}}{dy}\right]_0\right\} + \mathcal{O}(h^4)\end{aligned} \tag{48}$$

And $\left[d\bar{\phi}^X_{\langle i\pm\frac{1}{2}\rangle}/dy\right]_0$ may be approximated with central difference:

$$\left[\frac{d\bar{\phi}^X_{\langle i\pm\frac{1}{2}\rangle}}{dy}\right]_0 = \frac{\bar{\phi}^X_{\langle i\pm\frac{1}{2}\rangle}(y_n) - \bar{\phi}^X_{\langle i\pm\frac{1}{2}\rangle}(y_s)}{\Delta y} + \mathcal{O}(h^2) \tag{49}$$

The slice averages of $\bar{\phi}^X_{\langle i\pm\frac{1}{2}\rangle}$ on the right side of Equation (48) are computed again with the aid of the PKP theorem. Since

$$\begin{aligned}\bar{\phi}^{xy}_{\left[i\pm\frac{1}{2}\right]} = \frac{\int_{y_s}^{y_n}(\Delta x)_{\langle i\pm\frac{1}{2}\rangle}\bar{\phi}^X_{\langle i\pm\frac{1}{2}\rangle}dy}{(\Delta V)_{\left[i\pm\frac{1}{2}\right]}} &= \frac{1}{\overline{(\Delta x)}_{\left[i\pm\frac{1}{2}\right]}}\frac{\int_{y_s}^{y_n}(\Delta x)_{\langle i\pm\frac{1}{2}\rangle}\bar{\phi}^X_{\langle i\pm\frac{1}{2}\rangle}dy}{\Delta y} \\ &= \frac{1}{\overline{(\Delta x)}_{\left[i\pm\frac{1}{2}\right]}}\left\{\overline{(\Delta x)}_{\left[i\pm\frac{1}{2}\right]}\overline{\left(\bar{\phi}^X\right)}^y_{\left[i\pm\frac{1}{2}\right]} \right. \\ &\left. + \frac{\Delta y^2}{12}\left[\frac{d(\Delta x)_{\langle i\pm\frac{1}{2}\rangle}}{dy}\right]_0\left[\frac{d\bar{\phi}^X_{\langle i\pm\frac{1}{2}\rangle}}{dy}\right]_0\right\} + \mathcal{O}(h^3)\end{aligned} \tag{50}$$

it indicates

$$\overline{\left(\bar{\phi}^X\right)}^y_{\left[i\pm\frac{1}{2}\right]} = \bar{\phi}^{xy}_{\left[i\pm\frac{1}{2}\right]} - \frac{\Delta y^2}{12\overline{(\Delta x)}_{\left[i\pm\frac{1}{2}\right]}}\left[\frac{d(\Delta x)_{\langle i\pm\frac{1}{2}\rangle}}{dy}\right]_0\left[\frac{d\bar{\phi}^X_{\langle i\pm\frac{1}{2}\rangle}}{dy}\right]_0 + \mathcal{O}(h^3) \tag{51}$$

Assembling equations (45) through (51) together, we end up with

$$A_i\bar{\phi}^{y}_{(i-1)} + B_i\bar{\phi}^{y}_{(i)} + C_i\bar{\phi}^{y}_{(i+1)} = D_i\bar{\phi}^{xy}_{\left[i-\frac{1}{2}\right]} + E_i\bar{\phi}^{xy}_{\left[i+\frac{1}{2}\right]} + r_i + \mathcal{O}(h^3) \tag{52}$$

where

$$\begin{cases} A_i = 0 \\ B_i = \overline{(\Delta x)}_{\left[i-\frac{1}{2}\right]} + \overline{(\Delta x)}_{\left[i+\frac{1}{2}\right]} \\ C_i = 0 \\ D_i = \overline{(\Delta x)}_{\left[i+\frac{1}{2}\right]} \\ E_i = \overline{(\Delta x)}_{\left[i-\frac{1}{2}\right]} \\ r_i = \dfrac{\Delta y^2}{12}\Bigg\{\left(\left[\dfrac{d(\Delta x)_{\langle i+\frac{1}{2}\rangle}}{dy}\right]_0 - \dfrac{\overline{(\Delta x)}_{\left[i+\frac{1}{2}\right]}}{\overline{(\Delta x)}_{\left[i-\frac{1}{2}\right]}}\left[\dfrac{d(\Delta x)_{\langle i-\frac{1}{2}\rangle}}{dy}\right]_0\right)\left[\dfrac{d\bar{\phi}^{X}_{\langle i-\frac{1}{2}\rangle}}{dy}\right]_0 \\ + \left(\left[\dfrac{d(\Delta x)_{\langle i-\frac{1}{2}\rangle}}{dy}\right]_0 - \dfrac{\overline{(\Delta x)}_{\left[i-\frac{1}{2}\right]}}{\overline{(\Delta x)}_{\left[i+\frac{1}{2}\right]}}\left[\dfrac{d(\Delta x)_{\langle i+\frac{1}{2}\rangle}}{dy}\right]_0\right)\left[\dfrac{d\bar{\phi}^{X}_{\langle i+\frac{1}{2}\rangle}}{dy}\right]_0 - \left(\left[\dfrac{d(\Delta x)_{\langle i-\frac{1}{2}\rangle}}{dy}\right]_0 + \left[\dfrac{d(\Delta x)_{\langle i+\frac{1}{2}\rangle}}{dy}\right]_0\right)\left[\dfrac{d\phi_{(i)}}{dy}\right]_0\Bigg\} \end{cases} \tag{53}$$

This equation is valid for all interior vertical faces including those sandwiched between two rectangular cells, for which $r_i = 0$.

The development of equations for boundary faces are similar to the procedure presented above. And the resultant formulas are of the same general form as Equation (52). For a west (left) boundary face $(i)$ that only has neighboring cells to its east (right), we have

$$\begin{cases} A_i = 0 \\ B_i = 1 \\ C_i = 1 \\ D_i = 0 \\ E_i = 2 \\ r_i = -\dfrac{\Delta y^2}{6\overline{(\Delta x)}_{\left[i+\frac{1}{2}\right]}}\left[\dfrac{d(\Delta x)_{\langle i+\frac{1}{2}\rangle}}{dy}\right]_0\left[\dfrac{d\bar{\phi}^{X}_{\langle i+\frac{1}{2}\rangle}}{dy}\right]_0 \end{cases} \tag{54}$$

For a right boundary face $(i)$ that only has left neighboring cells, the coefficients are

$$\begin{cases} A_i = 1 \\ B_i = 1 \\ C_i = 0 \\ D_i = 2 \\ E_i = 0 \\ r_i = -\dfrac{\Delta y^2}{6\overline{(\Delta x)}_{\left[i-\frac{1}{2}\right]}}\left[\dfrac{d(\Delta x)_{\langle i-\frac{1}{2}\rangle}}{dy}\right]_0\left[\dfrac{d\bar{\phi}^{X}_{\langle i-\frac{1}{2}\rangle}}{dy}\right]_0 \end{cases} \tag{55}$$

Equations (52) through (55) and their counterparts in another coordinate direction can be used for faces along the same row or the same column of mesh, with the exception of boundary faces of solitary cells. For a solitary cell, the equations for its two boundary faces, Equations (54) and (55) are the same, and an extra equation must be provided, which will be discussed in the next section.

Sliding averages of derivatives can be related to cell averages in almost the same way elucidated above. For the slices shown in Figure 8, we have

$$\left(\frac{\partial\phi}{\partial x}\right)_i = f_i\phi_i + d_i\bar{\phi}^X_{\langle i-\frac{1}{2}\rangle} + e_i\bar{\phi}^X_{\langle i-\frac{1}{2}\rangle} + \mathcal{O}(h^2) \tag{56}$$

where

$$\begin{cases} f_i = \dfrac{2\left[(\Delta x)_{\langle i+\frac{1}{2}\rangle} - (\Delta x)_{\langle i-\frac{1}{2}\rangle}\right]}{(\Delta x)_{\langle i-\frac{1}{2}\rangle}(\Delta x)_{\langle i+\frac{1}{2}\rangle}} \\ d_i = -\dfrac{2(\Delta x)_{\langle i+\frac{1}{2}\rangle}}{(\Delta x)_{\langle i-\frac{1}{2}\rangle}\left[(\Delta x)_{\langle i-\frac{1}{2}\rangle} + (\Delta x)_{\langle i+\frac{1}{2}\rangle}\right]} \\ e_i = \dfrac{2(\Delta x)_{\langle i-\frac{1}{2}\rangle}}{(\Delta x)_{\langle i+\frac{1}{2}\rangle}\left[(\Delta x)_{\langle i-\frac{1}{2}\rangle} + (\Delta x)_{\langle i+\frac{1}{2}\rangle}\right]} \end{cases} \tag{57}$$

If $(i)$ is a west (left) boundary face, the formula is

$$b_i\left(\frac{\partial\phi}{\partial x}\right)_i + c_i\left(\frac{\partial\phi}{\partial x}\right)_{i+1} = -6\phi_i + 6\bar{\phi}^X_{\langle i+\frac{1}{2}\rangle} + \mathcal{O}(h^3) \tag{58}$$

where

$$\begin{cases} b_i = 2(\Delta x)_{\langle i+\frac{1}{2}\rangle} \\ c_i = (\Delta x)_{\langle i+\frac{1}{2}\rangle} \end{cases} \tag{59}$$

The corresponding formula for the east (right) boundary face $(i)$ is

$$a_i\left(\frac{\partial\phi}{\partial x}\right)_{i-1} + b_i\left(\frac{\partial\phi}{\partial x}\right)_i = 6\phi_i - 6\bar{\phi}^X_{\langle i-\frac{1}{2}\rangle} + \mathcal{O}(h^3) \tag{60}$$

where

$$\begin{cases} a_i = (\Delta x)_{\langle i-\frac{1}{2}\rangle} \\ b_i = 2(\Delta x)_{\langle i-\frac{1}{2}\rangle} \end{cases} \tag{61}$$

If the current face $(i)$ has an irregular west (left) neighboring cell and a rectangular east (right) neighboring cell, $(\Delta x)_{\langle i-\frac{1}{2}\rangle}$ vanishes at one end of the interface, see Figure 8. For such situations we can multiply both sides of Equation (56) by $(\Delta x)_{\langle i-\frac{1}{2}\rangle}$ then no singularities arise when we integrate the equation over the $y$-range of the face, which in virtue of the PKP theorem becomes

$$A_i \overline{\left(\frac{\partial \phi}{\partial x}\right)}^{y}_{(i-1)} + B_i \overline{\left(\frac{\partial \phi}{\partial x}\right)}^{y}_{(i)} + C_i \overline{\left(\frac{\partial \phi}{\partial x}\right)}^{y}_{(i+1)} = F_i \bar{\phi}^{y}_{(i)} + D_i \overline{(\phi)}^{xy}_{\left[i-\frac{1}{2}\right]} + E_i \overline{(\phi)}^{xy}_{\left[i+\frac{1}{2}\right]} + r_i + \mathcal{O}(h^3) \tag{62}$$

where

$$\left\{\begin{array}{c}
A_i = 0 \\
B_i = \overline{(\Delta x)}_{\left[i-\frac{1}{2}\right]} \\
C_i = 0 \\
F_i = 2\left[1 - \frac{\overline{(\Delta x)}_{\left[i-\frac{1}{2}\right]}}{(\Delta x)_{\left[i+\frac{1}{2}\right]}}\right] \\
D_i = -2(\Delta x)_{\left[i+\frac{1}{2}\right]} \overline{\left[\frac{1}{(\Delta x)_{\langle i-\frac{1}{2}\rangle} + (\Delta x)_{\left[i+\frac{1}{2}\right]}}\right]}^{y} \\
E_i = -D_i - F_i \\
r_i = \frac{\Delta y^2}{12}\left[\frac{d(\Delta x)_{\langle i-\frac{1}{2}\rangle}}{dy}\right]_0 \left\{-\frac{2}{(\Delta x)_{\left[i+\frac{1}{2}\right]}}\left[\frac{d\phi_{(i)}}{dy}\right]_0 + \left[\frac{2(\Delta x)_{\left[i+\frac{1}{2}\right]}}{\left\{(\Delta x)_{\langle i-\frac{1}{2}\rangle} + (\Delta x)_{\left[i+\frac{1}{2}\right]}\right\}^2_0} - \frac{D_i}{B_i}\right]\left[\frac{d\bar{\phi}^{X}_{\langle i-\frac{1}{2}\rangle}}{dy}\right]_0 \right. \\
\left. + \frac{2}{(\Delta x)_{\left[i+\frac{1}{2}\right]}}\left[\frac{(\Delta x)_{\langle i-\frac{1}{2}\rangle}\left\{(\Delta x)_{\langle i-\frac{1}{2}\rangle} + 2(\Delta x)_{\left[i+\frac{1}{2}\right]}\right\}}{\left\{(\Delta x)_{\langle i-\frac{1}{2}\rangle} + (\Delta x)_{\left[i+\frac{1}{2}\right]}\right\}^2}\right]_0 \left[\frac{d\bar{\phi}^{X}_{\langle i+\frac{1}{2}\rangle}}{dy}\right]_0 - \left[\frac{d\left(\frac{\partial \phi}{\partial x}\right)_{(i)}}{dy}\right]_0 \right\}
\end{array}\right. \tag{63}$$

Notice that $A_i$, $B_i$, $C_i$, $D_i$, $E_i$, $F_i$ and coefficients in $r_i$ are all constants and need to be calculated only once.

One may wonder why we do not multiply both sides of Equation (56) by $(\Delta x)_{\langle i-\frac{1}{2}\rangle}(\Delta x)_{\langle i+\frac{1}{2}\rangle}\left[(\Delta x)_{\langle i-\frac{1}{2}\rangle} + (\Delta x)_{\langle i+\frac{1}{2}\rangle}\right]$ then taking sliding average, so that the equation obtained would be universal for all interior vertical faces. This idea does not work because the equation then would have a fifth-order truncation error, which is not supported by the current nominally fourth-order, indeed third-order accurate version of PKP theorem.

For a face $(i)$ with an irregular east (right) neighboring boundary cell and a rectangular west (left) neighboring cell, we can multiply both sides of Equation (56) by $(\Delta x)_{\langle i+\frac{1}{2}\rangle}$ then integrating the equation over the $y$-range of the face. The coefficients and residual term of Equation (62) are

$$
\begin{cases}
A_i = 0 \\
B_i = \overline{(\Delta x)}_{\left[i+\frac{1}{2}\right]} \\
C_i = 0 \\
F_i = 2\left[\dfrac{\overline{(\Delta x)}_{\left[i+\frac{1}{2}\right]}}{(\Delta x)_{\left[i-\frac{1}{2}\right]}} - 1\right] \\
D_i = -E_i - F_i \\
E_i = 2(\Delta x)_{\left[i-\frac{1}{2}\right]} \overline{\left[\dfrac{1}{(\Delta x)_{\langle i+\frac{1}{2}\rangle} + (\Delta x)_{\left[i-\frac{1}{2}\right]}}\right]}^{y} \\
r_i = \dfrac{\Delta y^2}{12}\left[\dfrac{d(\Delta x)_{\langle i+\frac{1}{2}\rangle}}{dy}\right]_0 \left\{\dfrac{2}{(\Delta x)_{\left[i-\frac{1}{2}\right]}}\left[\dfrac{d\phi_{(i)}}{dy}\right]_0 - \left[\dfrac{2(\Delta x)_{\left[i-\frac{1}{2}\right]}}{\left\{(\Delta x)_{\langle i+\frac{1}{2}\rangle} + (\Delta x)_{\left[i-\frac{1}{2}\right]}\right\}_0^2} + \dfrac{E_i}{B_i}\right]\left[\dfrac{d\bar{\phi}^X_{\langle i+\frac{1}{2}\rangle}}{dy}\right]_0 \right. \\
\left. \qquad - \dfrac{2}{(\Delta x)_{\left[i-\frac{1}{2}\right]}}\left[\dfrac{(\Delta x)_{\langle i+\frac{1}{2}\rangle}\left\{(\Delta x)_{\langle i+\frac{1}{2}\rangle} + 2(\Delta x)_{\left[i-\frac{1}{2}\right]}\right\}}{\left\{(\Delta x)_{\langle i+\frac{1}{2}\rangle} + (\Delta x)_{\left[i-\frac{1}{2}\right]}\right\}^2}\right]_0 \left[\dfrac{d\bar{\phi}^X_{\langle i-\frac{1}{2}\rangle}}{dy}\right]_0 - \left[\dfrac{d\left(\frac{\partial\phi}{\partial x}\right)_{(i)}}{dy}\right]_0 \right\}
\end{cases}
\tag{64}
$$

Equations (63) and (64) (and their counterparts for horizontal faces) apply to all interfaces except those between twin cells because both $(\Delta x)_{\langle i+\frac{1}{2}\rangle}$ and $(\Delta x)_{\langle i-\frac{1}{2}\rangle}$ vanish at one end of such faces. This situation will be discussed shortly.

For a west (left) boundary face $(i)$, the coefficients and residual term of Equation (62) are

$$
\begin{cases}
A_i = 0 \\
B_i = 2\overline{(\Delta x)}_{\left[i+\frac{1}{2}\right]} \\
C_i = \overline{(\Delta x)}_{\left[i+\frac{1}{2}\right]} \\
F_i = -6 \\
D_i = 0 \\
E_i = 6 \\
r_i = -\dfrac{\Delta y^2}{12}\left[\dfrac{d(\Delta x)_{\langle i+\frac{1}{2}\rangle}}{dy}\right]_0 \left\{\dfrac{6}{\overline{(\Delta x)}_{\left[i+\frac{1}{2}\right]}}\left[\dfrac{d\bar{\phi}^X_{\langle i+\frac{1}{2}\rangle}}{dy}\right]_0 + 2\left[\dfrac{d\left(\frac{\partial\phi}{\partial x}\right)_{(i)}}{dy}\right]_0 + \left[\dfrac{d\left(\frac{\partial\phi}{\partial x}\right)_{(i+1)}}{dy}\right]_0\right\}
\end{cases}
\tag{65}
$$

For an east (right) boundary face $(i)$, the coefficients and residual term of Equation (62) are

$$
\begin{cases}
A_i = \overline{(\Delta x)}_{\left[i-\frac{1}{2}\right]} \\
B_i = 2\overline{(\Delta x)}_{\left[i-\frac{1}{2}\right]} \\
C_i = 0 \\
F_i = 6 \\
D_i = -6 \\
E_i = 0 \\
r_i = -\frac{\Delta y^2}{12}\left[\frac{d(\Delta x)_{\langle i-\frac{1}{2}\rangle}}{dy}\right]_0 \left\{-\frac{6}{\overline{(\Delta x)}_{\left[i-\frac{1}{2}\right]}}\left[\frac{d\bar{\phi}^X_{\langle i-\frac{1}{2}\rangle}}{dy}\right]_0 + 2\left[\frac{d\left(\frac{\partial\phi}{\partial x}\right)_{(i)}}{dy}\right]_0 + \left[\frac{d\left(\frac{\partial\phi}{\partial x}\right)_{(i-1)}}{dy}\right]_0\right\}
\end{cases}
\tag{66}
$$

Although Equations (65) and (66) are applicable to boundary faces of twin cells, we do not use them at such faces. Instead, we apply them to a twin-cell interface, say face $(i)$ and then add them up to construct an equation for the twin-cell interface, to which Equations (63) and (64) do not apply. The equation again assumes the form of Equation (62) and

$$
\begin{cases}
A_i = \overline{(\Delta x)}_{\left[i-\frac{1}{2}\right]} \\
B_i = 2\left[\overline{(\Delta x)}_{\left[i-\frac{1}{2}\right]} + \overline{(\Delta x)}_{\left[i+\frac{1}{2}\right]}\right] \\
C_i = \overline{(\Delta x)}_{\left[i+\frac{1}{2}\right]} \\
F_i = 0 \\
D_i = -6 \\
E_i = 6 \\
r_i = -\frac{\Delta y^2}{12}\left[\frac{d(\Delta x)_{\langle i-\frac{1}{2}\rangle}}{dy}\right]_0 \left\{-\frac{6}{\overline{(\Delta x)}_{\left[i-\frac{1}{2}\right]}}\left[\frac{d\bar{\phi}^X_{\langle i-\frac{1}{2}\rangle}}{dy}\right]_0 + 2\left[\frac{d\left(\frac{\partial\phi}{\partial x}\right)_{(i)}}{dy}\right]_0 + \left[\frac{d\left(\frac{\partial\phi}{\partial x}\right)_{(i-1)}}{dy}\right]_0\right\} \\
\quad -\frac{\Delta y^2}{12}\left[\frac{d(\Delta x)_{\langle i+\frac{1}{2}\rangle}}{dy}\right]_0 \left\{\frac{6}{\overline{(\Delta x)}_{\left[i+\frac{1}{2}\right]}}\left[\frac{d\bar{\phi}^X_{\langle i+\frac{1}{2}\rangle}}{dy}\right]_0 + 2\left[\frac{d\left(\frac{\partial\phi}{\partial x}\right)_{(i)}}{dy}\right]_0 + \left[\frac{d\left(\frac{\partial\phi}{\partial x}\right)_{(i+1)}}{dy}\right]_0\right\}
\end{cases}
\tag{67}
$$

Then we need equations for boundary faces of twin cells different from Equations (65) and (66). This issue is discussed in the next section.

## 2.6. Special treatments of solitary and twin cells

As mentioned in the previous section, only one of Equations (54) and (55) (or their counterparts in another direction) can be applied to a solitary cell when we calculate sliding averages along its boundaries, for example, when we calculate $\bar{\phi}^x$ along boundary faces 1 and 2 of solitary cell A in Figure 9. There are many ways, for instance polynomial fitting, to fix this problem. Hereby we present a method which is more aligned with the spirit of the current work.

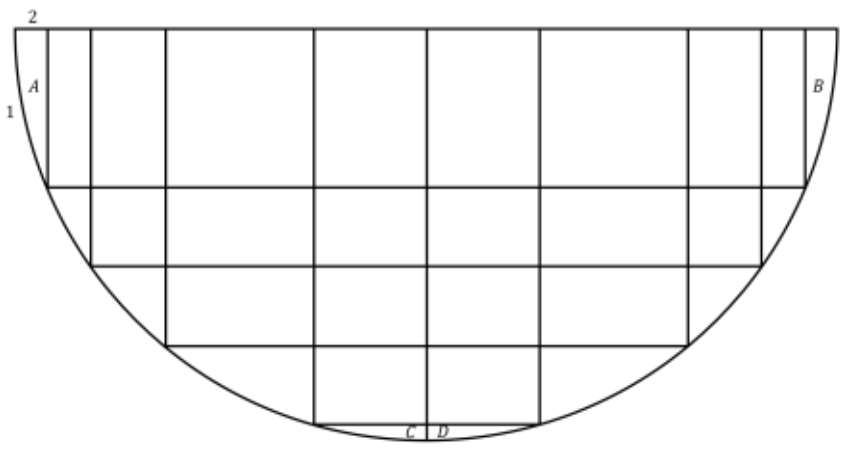

Figure 9

First, a solitary cell is only solitary in one direction. There is no difficulty in computing their sliding averages in the other direction. For example, we can find $\bar{\phi}^y$ along face 1 of solitary cell A in Figure 9 because this cell has enough neighboring cells along the $x$-direction, i.e., cell A is not solitary in the $x$-direction. Sliding average $\bar{\phi}^x$ along the same face may be evaluated from $\bar{\phi}^y$ as follows.

$$\begin{aligned}
\bar{\phi}^x_{(1)} &= \frac{\int_{x_{1,w}}^{x_{1,e}} \phi_1 dx_1}{(\Delta x)_{(1)}} = \frac{1}{(\Delta x)_{(1)}} \int_{y_{1,n}}^{y_{1,s}} \phi_1 \frac{dx_1}{dy_1} dy_1 = -\frac{1}{(\Delta x)_{(1)}} \int_{y_{1,s}}^{y_{1,n}} \phi_1 \frac{dx_1}{dy_1} dy_1 \\
&= -\frac{(\Delta y)_{(1)}}{(\Delta x)_{(1)}} \overline{\left(\phi_1 \frac{dx_1}{dy_1}\right)}^y \\
&= -\frac{(\Delta y)_{(1)}}{(\Delta x)_{(1)}} \left[ \overline{\left(\frac{dx_1}{dy_1}\right)}^y \bar{\phi}^y_{(1)} + \frac{(\Delta y)^2_{(1)}}{12} \left(\frac{d^2 x_1}{dy_1^2}\right)_0 \left[\frac{d\phi_{(1)}}{dy}\right]_0 \right. \\
&\quad \left. + \mathcal{O}(h^4) \right]
\end{aligned} \tag{68}$$

in which

$$\overline{\left(\frac{dx_1}{dy_1}\right)}^y = \frac{1}{(\Delta y)_{(1)}} \int_{y_{1,s}}^{y_{1,n}} \frac{dx_1}{dy_1} dy_1 = \frac{1}{(\Delta y)_{(1)}} \int_{x_{1,e}}^{x_{1,w}} dx_1 = -\frac{(\Delta x)_{(1)}}{(\Delta y)_{(1)}} \tag{69}$$

Therefore,

$$\bar{\phi}^x_{(1)} = \bar{\phi}^y_{(1)} - \frac{(\Delta y)_{(1)}}{(\Delta x)_{(1)}} \left\{ \frac{(\Delta y)^2_{(1)}}{12} \left(\frac{d^2 x_1}{dy_1^2}\right)_0 \left[\frac{d\phi_{(1)}}{dy}\right]_0 + \mathcal{O}(h^4) \right\} \tag{70}$$

and $(d^2 x_1/dy_1^2)_0$ can be computed exactly and $\left[d\phi_{(1)}/dy\right]_0$ may be estimated with methods introduced in Section 2.4. Once we find $\bar{\phi}^x_{(1)}$, the slice average along the other boundary face, face 2, may be evaluated by using an equation like (54) or (55).

Although unnecessary, the sliding average of derivative $\overline{(\partial\phi/\partial y)}^x$ along face 1 also can be evaluated from $\overline{(\partial\phi/\partial x)}^y$ with the help of Equation (27) as follows:

$$
\begin{aligned}
\overline{\left(\frac{\partial\phi}{\partial y}\right)}^{x}_{(1)} &= \frac{1}{(\Delta x)_{(1)}}\int_{x_{1,w}}^{x_{1,e}}\left(\frac{\partial\phi}{\partial y}\right)_1 dx_1 \\
&= \frac{1}{(\Delta x)_{(1)}}\int_{y_{1,n}}^{y_{1,s}}\left[\frac{d\phi_1}{dy_1}-\left(\frac{\partial\phi}{\partial x}\right)_1\frac{dx_1}{dy_1}\right]\frac{dx_1}{dy_1}dy_1 \\
&= \frac{1}{(\Delta x)_{(1)}}\int_{y_{1,s}}^{y_{1,n}}\left[\left(\frac{\partial\phi}{\partial x}\right)_1\frac{dx_1}{dy_1}-\frac{d\phi_1}{dy_1}\right]\frac{dx_1}{dy_1}dy_1 \\
&= \frac{1}{(\Delta x)_{(1)}}\left\{\overline{\left[\left(\frac{\partial\phi}{\partial x}\right)_1\left(\frac{dx_1}{dy_1}\right)^2\right]}^{y}(\Delta y)_{(1)}-\left(\phi_1\frac{dx_1}{dy_1}\right)\Bigg|_{y_{1,s}}^{y_{1,n}}\right. \\
&\left.+\overline{\left(\phi_1\frac{d^2x_1}{dy_1^2}\right)}^{y}(\Delta y)_{(1)}\right\}
\end{aligned}
\tag{71}
$$

The sliding averages of products are then evaluated using the PKP theorem, which gives

$$
\begin{aligned}
\overline{\left(\frac{\partial\phi}{\partial y}\right)}^{x}_{(1)} = \frac{(\Delta y)_{(1)}}{(\Delta x)_{(1)}}&\left\{\overline{\left(\frac{\partial\phi}{\partial x}\right)}^{y}_{(1)}\overline{\left[\left(\frac{dx_1}{dy_1}\right)^2\right]}^{y}+\frac{(\Delta y)^2_{(1)}}{6}\left[\frac{d\left(\frac{\partial\phi}{\partial x}\right)_{(1)}}{dy}\right]_0\left(\frac{dx_1}{dy_1}\frac{d^2x_1}{dy_1^2}\right)_0\right. \\
&-\frac{1}{(\Delta y)_{(1)}}\left(\phi_1\frac{dx_1}{dy_1}\right)\Bigg|_{y_{1,s}}^{y_{1,n}}+\frac{1}{(\Delta y)_{(1)}}\bar{\phi}^{y}_{(1)}\left(\frac{dx_1}{dy_1}\right)\Bigg|_{y_{1,s}}^{y_{1,n}} \\
&\left.+\frac{(\Delta y)^2_{(1)}}{12}\left[\frac{d\phi_{(1)}}{dy}\right]_0\left(\frac{d^3x_1}{dy_1^3}\right)_0+\mathcal{O}(h^4)\right\}
\end{aligned}
\tag{72}
$$

Terms on the right side of this equation are either calculated exactly or assessed with methods introduced in Sections 2.4 and 2.5. The same method is applied to boundary faces of twin cells, e.g., cells C and D in Figure 9.

## 2.7. Solving process

Momentum Equations (21) and (22), once sliding averages in them are represented by corresponding cell averages, may be written in the symbolic form of

$$
\frac{\partial}{\partial t}\bar{u}^{xy}_{[i]}+A^u_{[i]}\bar{u}^{xy}_{[i]}=\sum_{l|i}A^u_{[l]}\bar{u}^{xy}_{[l]}+r^u_{[i]}-\overline{\left(\frac{\partial p}{\partial x}\right)}^{xy}_{[i]}
\tag{73}
$$

and

$$
\frac{\partial}{\partial t}\bar{v}^{xy}_{[j]}+A^v_{[j]}\bar{v}^{xy}_{[j]}=\sum_{k|j}A^v_{[k]}\bar{v}^{xy}_{[k]}+r^v_{[j]}-\overline{\left(\frac{\partial p}{\partial y}\right)}^{xy}_{[j]}
\tag{74}
$$

These equations and the continuity equation can be solved by any collocated mesh Navier-Stokes solvers. Our method is as follows.

Since temporal accuracy is not the concern of the current study, we adopt a simple implicit scheme for time advancing. Equations (73) and (74) become

$$\frac{\left(\bar{u}^{xy}_{[i]}\right)^* - \left(\bar{u}^{xy}_{[i]}\right)^n}{\Delta t} + A^u_{[i]}\left(\bar{u}^{xy}_{[i]}\right)^* = \sum_{l|i} A^u_{[l]}\left(\bar{u}^{xy}_{[l]}\right)^* + \left(r^u_{[i]}\right)^m - \left[\overline{\left(\frac{\partial p}{\partial x}\right)}^{xy}_{[i]}\right]^m \tag{75}$$

and

$$\frac{\left(\bar{v}^{xy}_{[j]}\right)^* - \left(\bar{v}^{xy}_{[j]}\right)^n}{\Delta t} + A^v_{[j]}\left(\bar{v}^{xy}_{[j]}\right)^* = \sum_{k|j} A^v_{[k]}\left(\bar{v}^{xy}_{[k]}\right)^* + \left(r^v_{[j]}\right)^m - \left[\overline{\left(\frac{\partial p}{\partial y}\right)}^{xy}_{[j]}\right]^m \tag{76}$$

where the superscript $*$ indicates an intermediate flow field; $n$ denotes the previous time level and $m$ represents the previous iteration. To avoid odd-even decoupling of pressure and velocity, we develop a solving procedure inspired by the work of (Pascau, 2011). First, we define

$$\left(\bar{u}^{xy}_{[i]}\right)^{**} = \frac{\sum_{l|i} A^u_{[l]}\left(\bar{u}^{xy}_{[l]}\right)^* + \left(r^u_{[i]}\right)^m}{A^u_{[i]}}, \quad \left(\bar{v}^{xy}_{[j]}\right)^{**} = \frac{\sum_{k|j} A^v_{[k]}\left(\bar{v}^{xy}_{[k]}\right)^* + \left(r^v_{[j]}\right)^m}{A^v_{[j]}} \tag{77}$$

When steady state is reached, we have

$$\left(\bar{u}^{xy}_{[i]}\right)^* = \left(\bar{u}^{xy}_{[i]}\right)^{**} - \frac{1}{A^u_{[i]}}\overline{\left(\frac{\partial p}{\partial x}\right)}^{xy}_{[i]}, \quad \left(\bar{v}^{xy}_{[j]}\right)^* = \left(\bar{v}^{xy}_{[j]}\right)^{**} - \frac{1}{A^v_{[j]}}\overline{\left(\frac{\partial p}{\partial y}\right)}^{xy}_{[j]} \tag{78}$$

We then change Equation (52) a bit to calculate the advective sliding-averaged $u$-velocity along vertical interior faces if the steady-state is achieved:

$$\left[\bar{u}^y_{(i)}\right]^* = \frac{D_i}{B_i}\left(\bar{u}^{xy}_{\left[i-\frac{1}{2}\right]}\right)^{**} + \frac{E_i}{B_i}\left(\bar{u}^{xy}_{\left[i+\frac{1}{2}\right]}\right)^{**} - \frac{1}{B_i}\left(\frac{D_i}{A^u_{\left[i-\frac{1}{2}\right]}} + \frac{E_i}{A^u_{\left[i+\frac{1}{2}\right]}}\right)\overline{\left(\frac{\partial p}{\partial x}\right)}^{y}_{(i)} + \frac{r_i}{B_i} + \mathcal{O}(h^2) \tag{79}$$

The constant coefficients $B_i$, $D_i$, $E_i$ and the form of the residual term $r_i$ are not changed, see Equation (53). The change we made is replacing the pressure gradients of cells with the pressure gradient at the face, which is the basic idea of the Rhie-Chow momentum interpolation method (Rhie & Chow, 1983) which amounts to adding a third-order term to damp the possible pressure oscillations. Comparison of Equations (78) and (79) implies that

$$\frac{1}{A^u_{(i)}} = \frac{1}{B_i}\left(\frac{D_i}{A^u_{\left[i-\frac{1}{2}\right]}} + \frac{E_i}{A^u_{\left[i+\frac{1}{2}\right]}}\right) \tag{80}$$

We then construct a fictitious equation for $\left[\bar{u}^y_{(i)}\right]^*$ before the steady state arrives:

$$\frac{\left[\bar{u}^y_{(i)}\right]^* - \left[\bar{u}^y_{(i)}\right]^n}{\Delta t} + A^u_{(i)}\left[\bar{u}^y_{(i)}\right]^* = \sum_{l|i} A^u_{(l)}\left[\bar{u}^y_{(l)}\right]^* + \left[r^u_{(i)}\right]^m - \left[\overline{\left(\frac{\partial p}{\partial x}\right)}^{y}_{(i)}\right]^{m+1} \tag{81}$$

where

$$\sum_{l|i} A^u_{(l)}\left[\bar{u}^y_{(l)}\right]^* = A^u_{(i)}\left[\frac{D_i}{B_i}\left(\bar{u}^{xy}_{\left[i-\frac{1}{2}\right]}\right)^{**} + \frac{E_i}{B_i}\left(\bar{u}^{xy}_{\left[i+\frac{1}{2}\right]}\right)^{**}\right], \quad \left[r^u_{(i)}\right]^m = A^u_{(i)}\frac{(r_i)^m}{B_i} \tag{82}$$

A similar equation exists for $\left[\bar{v}_{(j)}^{x}\right]^{*}$ along horizontal interior faces:

$$\frac{\left[\bar{v}_{(j)}^{x}\right]^{*}-\left[\bar{v}_{(j)}^{x}\right]^{n}}{\Delta t}+A_{(j)}^{v}\left[\bar{v}_{(j)}^{x}\right]^{*}=\sum_{k|j}A_{(k)}^{v}\left[\bar{v}_{(k)}^{x}\right]^{*}+\left[r_{(j)}^{v}\right]^{m}-\left[\overline{\left(\frac{\partial p}{\partial y}\right)}_{(j)}^{x}\right]^{m+1} \tag{83}$$

These two equations can be written in more compact forms:

$$\left[\bar{u}_{(i)}^{y}\right]^{*}=-\frac{1}{\tilde{A}_{(i)}^{u}}\left[\overline{\left(\frac{\partial p}{\partial x}\right)}_{(i)}^{y}\right]^{m+1}+\left[R_{(i)}^{u}\right]^{m},\quad \left[\bar{v}_{(j)}^{x}\right]^{*}=-\frac{1}{\tilde{A}_{(j)}^{v}}\left[\overline{\left(\frac{\partial p}{\partial y}\right)}_{(j)}^{x}\right]^{m+1}+\left[R_{(j)}^{v}\right]^{m} \tag{84}$$

in which

$$\tilde{A}_{(i)}^{u}=A_{(i)}^{u}+\frac{1}{\Delta t},\quad \left[R_{(i)}^{u}\right]^{m}=\frac{1}{\tilde{A}_{(i)}^{u}}\left(\sum_{l|i}A_{(l)}^{u}\left[\bar{u}_{(l)}^{y}\right]^{*}+\left[r_{(i)}^{u}\right]^{m}+\frac{1}{\Delta t}\left[\bar{u}_{(i)}^{y}\right]^{n}\right) \tag{85}$$

and

$$\tilde{A}_{(j)}^{v}=A_{(j)}^{v}+\frac{1}{\Delta t},\quad \left[R_{(j)}^{v}\right]^{m}=\frac{1}{\tilde{A}_{(j)}^{v}}\left(\sum_{k|j}A_{(k)}^{v}\left[\bar{v}_{(k)}^{x}\right]^{*}+\left[r_{(j)}^{v}\right]^{m}+\frac{1}{\Delta t}\left[\bar{v}_{(j)}^{x}\right]^{n}\right) \tag{86}$$

Then we enforce the continuity equation, Equation (3) at each cell. Integration of the continuity equation over a cell, say cell $\left[i-\frac{1}{2},j-\frac{1}{2}\right]$ results in

$$\frac{\left[\bar{u}_{\left(i,j-\frac{1}{2}\right)}^{y}\right]^{*}-\left[\bar{u}_{\left(i-1,j-\frac{1}{2}\right)}^{y}\right]^{*}}{\overline{\Delta x}}+\frac{\left[\bar{v}_{\left(i-\frac{1}{2},j\right)}^{x}\right]^{*}-\left[\bar{v}_{\left(i-\frac{1}{2},j-1\right)}^{x}\right]^{*}}{\overline{\Delta y}}=0 \tag{87}$$

Substitution of Equations (81) and (83) into Equation (87) gives

$$\begin{aligned}
&\frac{1}{\overline{\overline{\Delta x}}}\left(\frac{1}{\tilde{A}_{\left(i,j-\frac{1}{2}\right)}^{u}}\left[\overline{\left(\frac{\partial p}{\partial x}\right)}_{\left(i,j-\frac{1}{2}\right)}^{y}\right]^{m+1}-\frac{1}{\tilde{A}_{\left(i-1,j-\frac{1}{2}\right)}^{u}}\left[\overline{\left(\frac{\partial p}{\partial x}\right)}_{\left(i-1,j-\frac{1}{2}\right)}^{y}\right]^{m+1}\right)\\
&\quad+\frac{1}{\overline{\overline{\Delta y}}}\left(\frac{1}{\tilde{A}_{\left(i-\frac{1}{2},j\right)}^{v}}\left[\overline{\left(\frac{\partial p}{\partial y}\right)}_{\left(i-\frac{1}{2},j\right)}^{x}\right]^{m+1}\right.\\
&\quad\left.-\frac{1}{\tilde{A}_{\left(i-\frac{1}{2},j-1\right)}^{v}}\left[\overline{\left(\frac{\partial p}{\partial y}\right)}_{\left(i-\frac{1}{2},j-1\right)}^{x}\right]^{m+1}\right)\\
&\quad=\frac{1}{\overline{\overline{\Delta x}}}\left(\left[R_{\left(i,j-\frac{1}{2}\right)}^{u}\right]^{m}-\left[R_{\left(i-1,j-\frac{1}{2}\right)}^{u}\right]^{m}\right)\\
&\quad+\frac{1}{\overline{\overline{\Delta y}}}\left(\left[R_{\left(i-\frac{1}{2},j\right)}^{v}\right]^{m}-\left[R_{\left(i-\frac{1}{2},j-1\right)}^{v}\right]^{m}\right)
\end{aligned} \tag{88}$$

which solves for the new cell-averaged pressure field $(\bar{p}^{xy})^{m+1}$ with the help of, say Equation (62). This new pressure field can be used to calculate a new intermediate velocity field. We may continue such iterations until a converged flow field, the flow field at $n+1$ time level, is reached. This time marching is repeated until the steady state solution is obtained.

# 3. Numerical Results

The lid-driven semi-circular cavity flow at Reynolds numbers $Re_D = 100$ and $Re_D = 1000$ are simulated with the current method. The Reynolds number is defined as

$$Re_D = \frac{UD}{\nu} \tag{89}$$

where $U$ is the lid speed, $D$ the cavity diameter and $\nu$ is the fluid kinematic viscosity. The mesh used in simulations is as shown in Figure 10. The upper half of the mesh is uniform in the $y$-direction; the lower half is uniform in the $x$-direction. The dividing line of these two halves is $y/D = -\cos(\pi/4)/2$, which is the bold horizontal line in the figure. The number (say $N$) of mesh cells along the lid ($y/D = 0$) is twice that along the vertical centerline ($x/D = 0$) of the cavity. In this way no matter what mesh densities are chosen, the lines $x/D = -\cos(\pi/4)/2$, and $y/D = -\cos(\pi/4)/2$ are always grid lines and the mesh along these lines are always almost uniform, so that results along these lines might be conveniently utilized to determine the order of accuracy of the numerical scheme. These two lines have at least one end at the circular boundary, which allows us to explore the effectiveness of the current complex boundary treatment method.

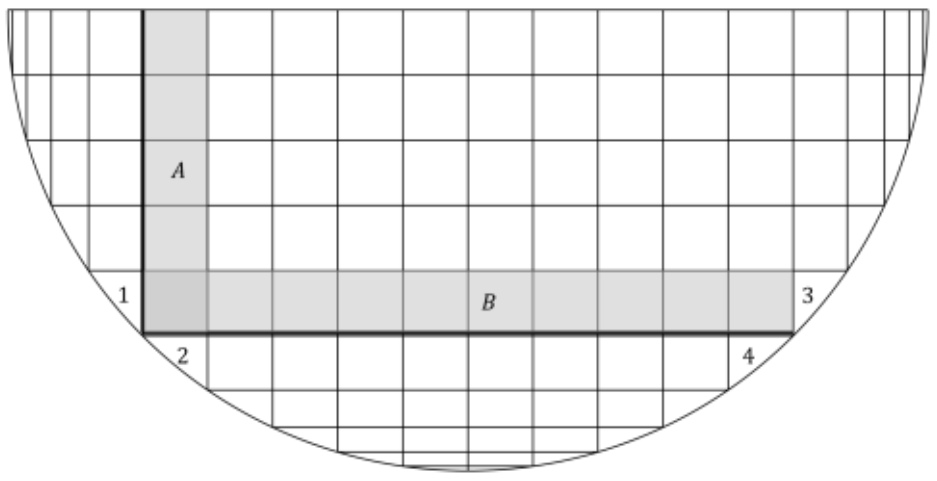

Figure 10

The common understanding of the order of accuracy is that if the discretization error at a fixed point in the flow field is proportional to the $n$th-power of the local cell size $h$, the numerical scheme is of $n$th-order accuracy. Although seemingly straightforward, this definition is not free of ambiguity. For example, what is $h$ of a rectangular cell? Is it $\Delta x$, or $\Delta y$, or $\sqrt{\Delta x \Delta y}$, or $(\Delta x + \Delta y)/2$, or something else? Such ambiguity reflects the fact that the definition of order of accuracy is indeed numerical scheme dependent. The definition of $h$ of the current method varies from formula to formula but is always of the form of $(\Delta x)^m(\Delta y)^n$ where $m + n = 1$. Therefore, for the current second-order accurate method the discretization error at a fixed point in the flow field is quartered when $\Delta x$ and $\Delta y$ of the local cell are both halved. However, it is impossible to halve $\Delta x$ and $\Delta y$ simultaneously for most irregular cells. For example, for the two twin cells, $\Delta y \approx \Delta x^2/D$ as $\Delta x \ll D/2$. Therefore if $\Delta x$ is halved, $\Delta y$ will be quartered for these cells, which may cause some apparent inconsistency in the order of accuracy. For instance, if the discretization error of a certain parameter at a twin cell is proportional to $\Delta x^4 \Delta y^{-2}$, which is still of second-order accuracy according to our definition, it may looks like a zeroth-order scheme because the discretization error is of the order of $\Delta x^4 \Delta y^{-2} \sim \Delta x^4 \Delta x^{-4} \sim \Delta x^0$. To avoid such subtleties, we use lines and cells shown in Figure 10 to facilitate the determination of order of accuracy as each of these cells has almost the same $\Delta x$ as $\Delta y$. Since a finite volume scheme usually results in a globally uniform order of accuracy (Kobayashi, 1999), the present practice should be able to reveal the actual order of accuracy of the current method.

Three mesh densities are tested for each Reynolds number. For $Re_D = 100$, they are $N = 20$, 40, and 60; for $Re_D = 1000$, they are $N = 40$, 60, and 80.

The face-averaged velocity distributions along the vertical line $x/D = -\cos(\pi/4)/2$ and the horizontal line $y/D = -\cos(\pi/4)/2$ at $Re_D = 100$ is shown in Figure 11. Symbols are current numerical results, and curves are benchmark data obtained by using commercial CFD software

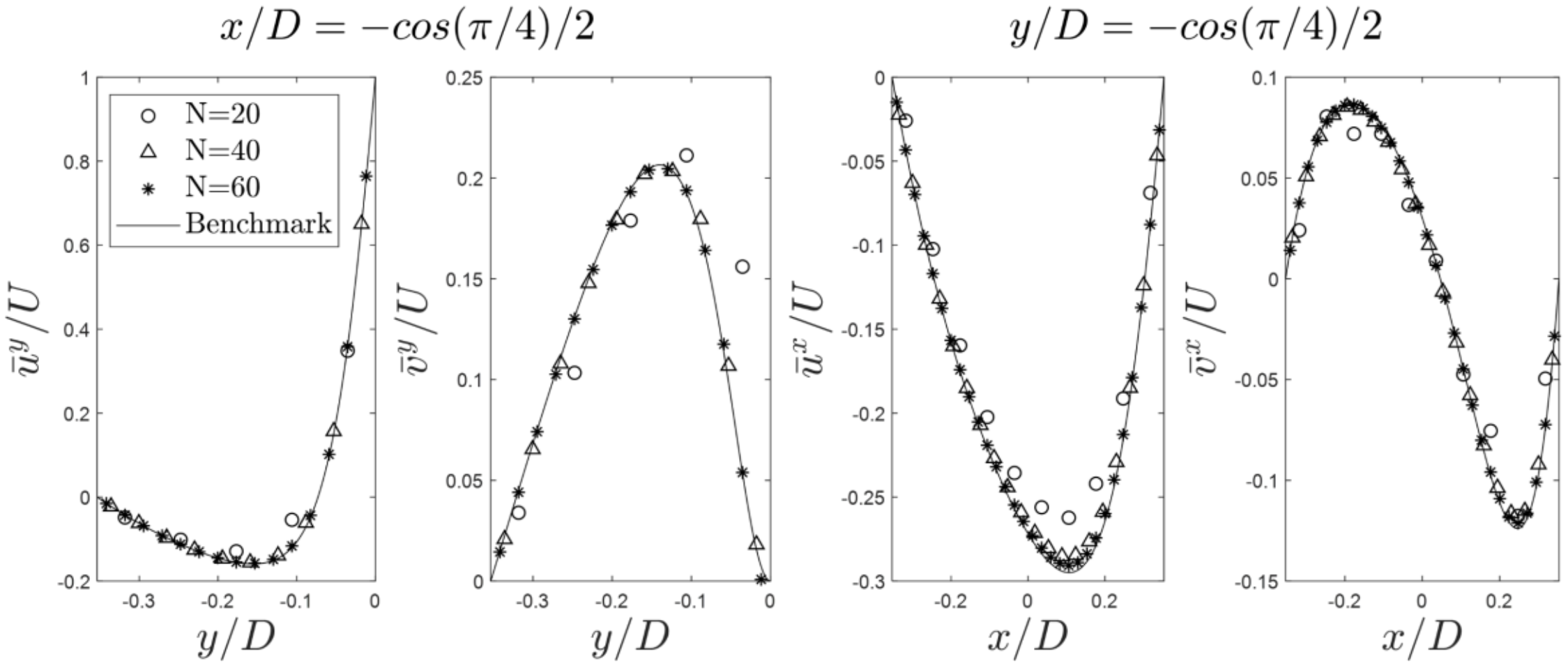

Figure 11

ANSYS FLUENT with a 43924-cell unstructured mesh. The face-averaged velocity gradients along these two lines are shown in Figure 12.

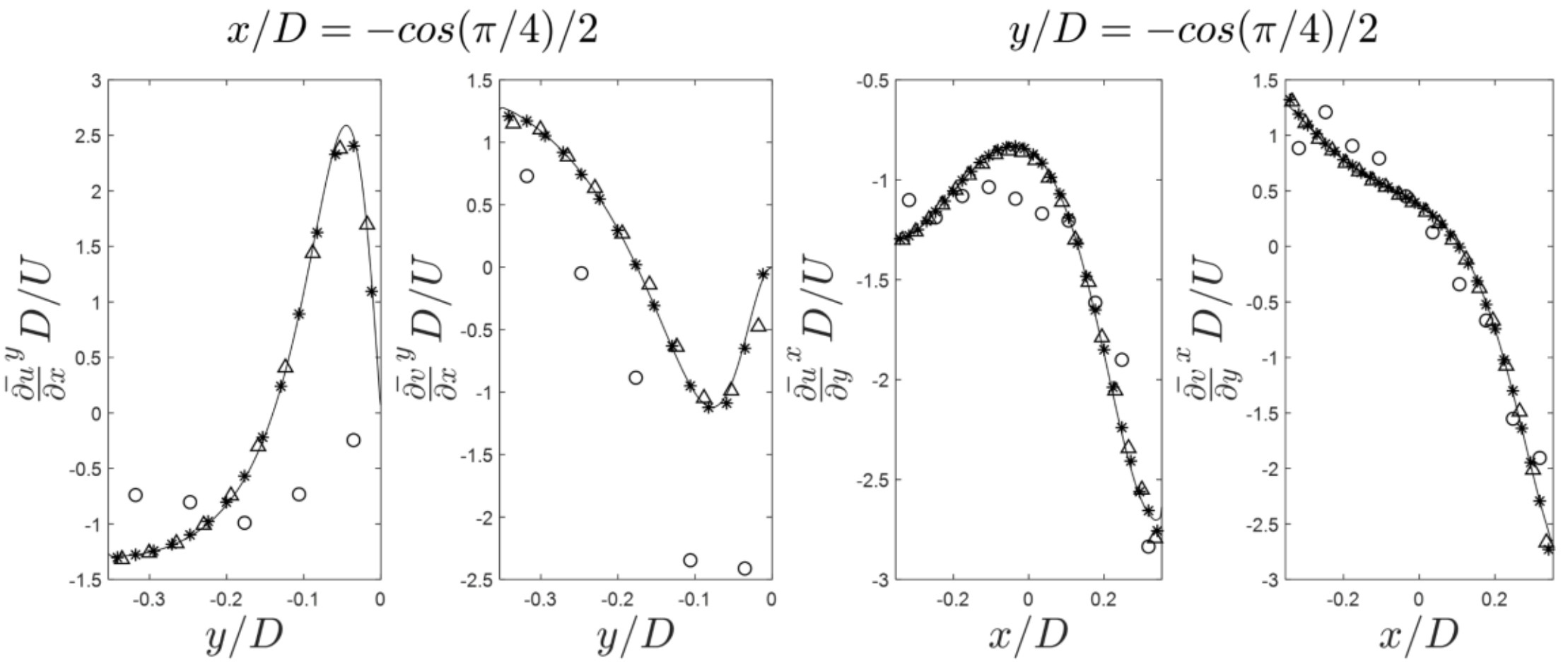

Figure 12

The cell-averaged velocities of the cells marked “$A$” in Figure 10 are compared with the corresponding benchmark data in Figure 13. Notice that the legends like “$N = 20$, Benchmark” in the figure denote benchmark data corresponding to the current numerical results with the specific mesh density.

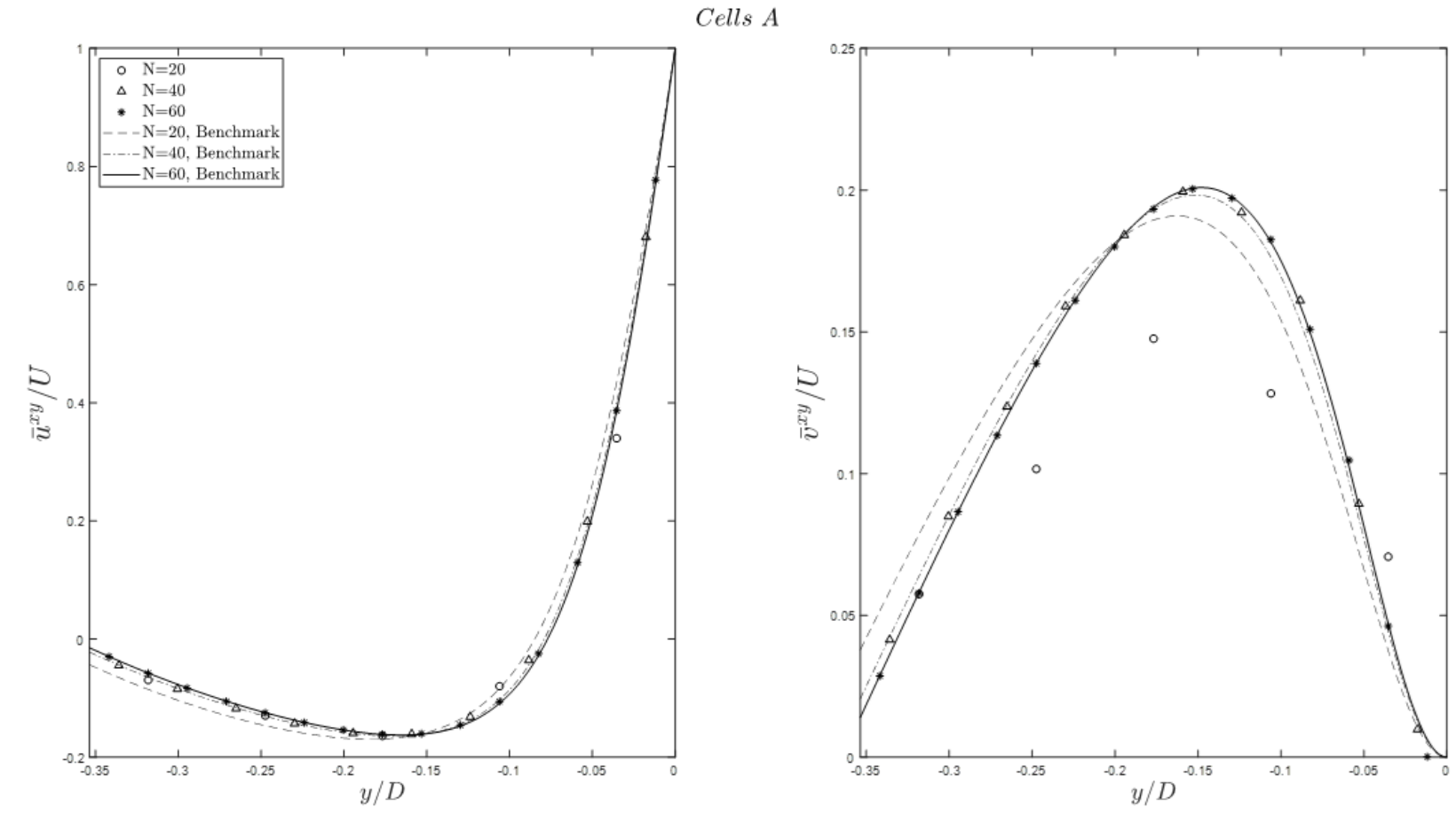

*Figure 13*

The cell-averaged velocities of the cells marked “$B$” in Figure 10 are compared with the corresponding benchmark data in Figure 14.

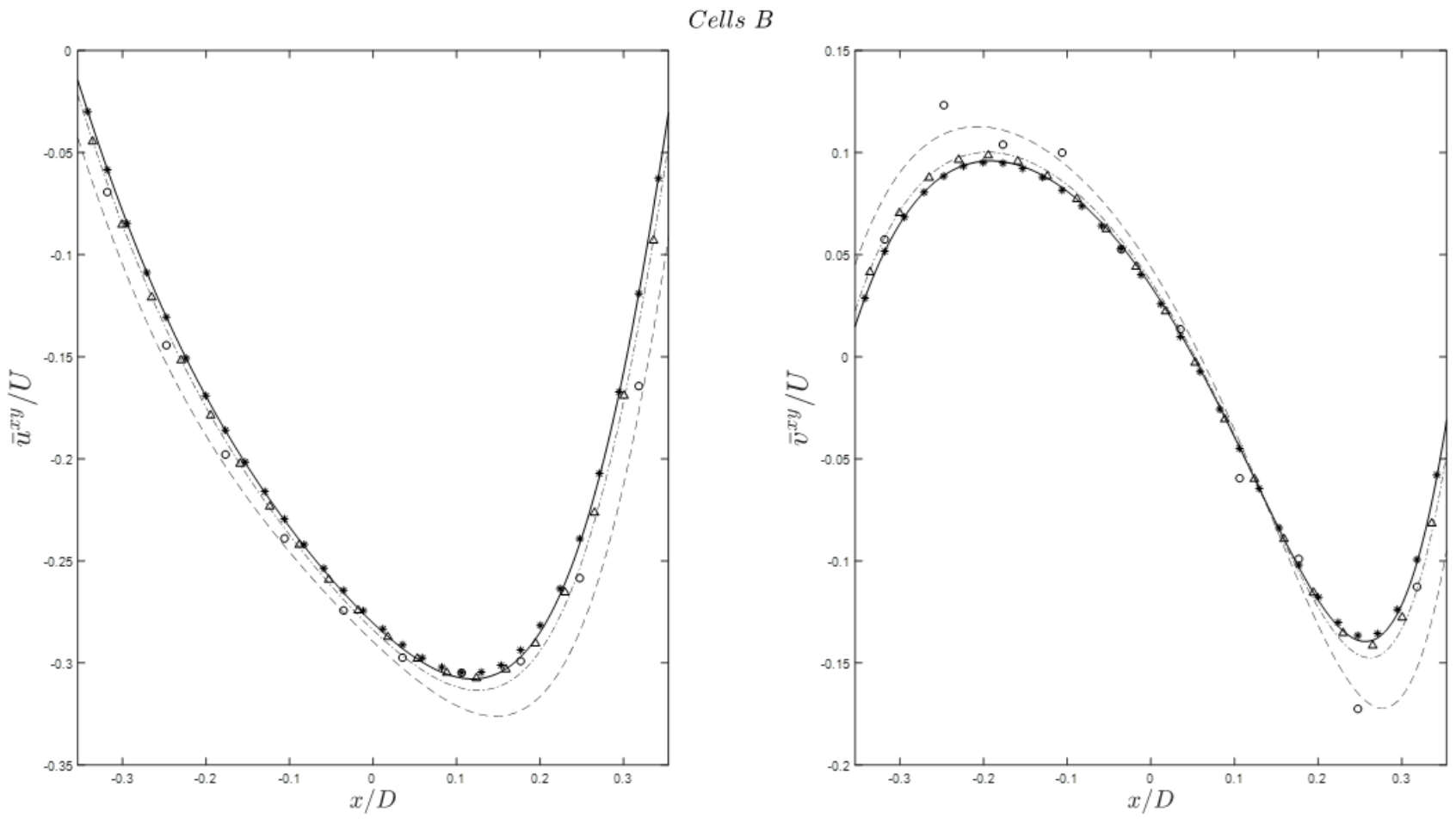

*Figure 14*

The cell-averaged velocity values of the four irregular boundary cells marked “1” through “4” are compared with the benchmark data in Table 1. The benchmark cell-averaged velocity values in the regions corresponding to cells 1 through 4 are obtained by using the zone separation and volume averaging ability of ANSYS FLUENT.

Table 1

| Velocity | Mesh ($N$) | Cell 1 | Cell 2 | Cell 3 | Cell 4 |
|---|---|---|---|---|---|
| $\bar{u}^{xy}/U$ | 20 | -0.0462 | -0.0124 | -0.0899 | -0.0375 |
| | Benchmark | -0.0278 | -0.0286 | -0.0611 | -0.0572 |
| | 40 | -0.0150 | -0.0154 | -0.0329 | -0.0315 |
| | Benchmark | -0.0146 | -0.0147 | -0.0310 | -0.0294 |
| | 60 | -0.0101 | -0.0101 | -0.0218 | -0.0211 |
| | Benchmark | -0.0098 | -0.0101 | -0.0207 | -0.0202 |
| $\bar{v}^{xy}/U$ | 20 | 0.0334 | 0.0126 | -0.0769 | -0.0273 |
| | Benchmark | 0.0323 | 0.0228 | -0.0695 | -0.0439 |
| | 40 | 0.0153 | 0.0131 | -0.0319 | -0.0251 |
| | Benchmark | 0.0158 | 0.0131 | -0.0331 | -0.0258 |
| | 60 | 0.0103 | 0.0090 | -0.0214 | -0.0180 |
| | Benchmark | 0.0104 | 0.0094 | -0.0216 | -0.0185 |

The flow field (streamlines) calculated by the current method with the finest mesh ($N = 60$) is compared with the benchmark result in Figure 15.

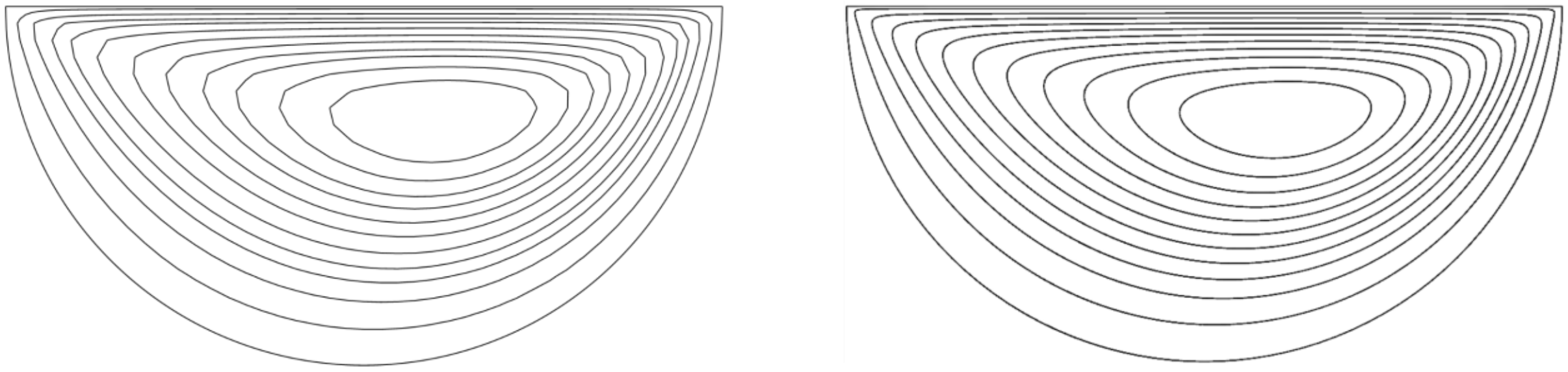

Figure 15

The arithmetic mean of discretization errors of face-averaged velocities and velocity gradients along lines $x/D = -\cos(\pi/4)/2$ and $y/D = -\cos(\pi/4)/2$, and those of cell-averaged velocities of cells "$A$" and "$B$", as well as mean errors of cell-averaged velocities of irregular boundary cells 1 through 4 shown in Figure 10 are used to evaluate the order of accuracy of the current method. The results are summarized in Table 2. All flow field variables show second or higher than second order error reduction as the mesh density is increased from $N = 20$ to $N = 40$. The majority of flow field variables show second or close to second-order behavior as the mesh is refined from $N = 40$ to $N = 60$. A few variables perform differently, most noticeably $\overline{\partial u/\partial x}^{y}$ along line $x/D = -\cos(\pi/4)/2$ which only has a 0.69th-order reduction of error as the mesh density is increased from $N = 40$ to $N = 60$, yet a staggering 4.77th-order of error deduction when the mesh density changes from $N = 20$ to $N = 40$. We do not have good explanations to this phenomenon. One may notice that all cell-averaged velocities, including those of irregular boundary cells show second or close to second-order accuracy.

Numerical results for the $Re_D = 1000$ case are compared with the benchmark in Figure 16 through Figure 20 and Table 3, with exactly the same order as those presented for the $Re_D = 100$ case. The mean errors of variables of interest are shown in Table 4. As can be observed, all flow variables show second or close to second-order (all higher than second-order, except one with a 1.99th-order) behavior as the mesh is refined from $N = 60$ to $N = 80$ with only one exception, namely $\bar{v}^y$ along line $x/D = -\cos(\pi/4)/2$, which shows a 1.84th-order, which is still very close to the second order accuracy. More importantly, all cell-averaged velocities, including those of irregular boundary cells reach second-order accuracy as the mesh density is increased from $N = 60$ to $N = 80$. It seems that for this high Reynolds number case the $N = 40$ mesh is too coarse, so that the leading term of the truncation error still does not dominate, and the second-order behavior does not manifest itself fully until the mesh density is high enough, say $N = 60$.

A few other parameters of interest are compared with the benchmark work of (Glowinski, Guidoboni, & Pan, 2006) in **Error! Reference source not found.**. These parameters include the value and location of the minimum stream function (the "eye" of the major vortex), and angles of detachment of the secondary vortex (see Figure 20 for definitions of these angles). Only parameters of the $Re_D = 1000$ case are presented here as the $Re_D = 100$ case was not discussed in the cited benchmark work. As can be observed, the current results which are obtained with the $N = 80$ mesh agree with the benchmark data very well.

*Table 2*

| Variable | Mesh ($N$) | Error | Order of accuracy | Variable | Mesh ($N$) | Error | Order of accuracy |
|---|---|---|---|---|---|---|---|
| $\bar{u}^y/U$ | 20 | 0.0239 | - | $\bar{v}^y/U$ | 20 | 0.0352 | - |
| | 40 | 0.0033 | 2.86 | | 40 | 0.0032 | 3.45 |
| | 60 | 0.0015 | 1.94 | | 60 | 0.0016 | 1.69 |
| $\bar{u}^x/U$ | 20 | 0.0228 | - | $\bar{v}^x/U$ | 20 | 0.0105 | - |
| | 40 | 0.0047 | 2.28 | | 40 | 0.0025 | 2.06 |
| | 60 | 0.0021 | 1.97 | | 60 | 0.0011 | 2.02 |
| $\overline{\frac{\partial u}{\partial x}}^y D/U$ | 20 | 1.1190 | - | $\overline{\frac{\partial v}{\partial x}}^y D/U$ | 20 | 1.0736 | - |
| | 40 | 0.0409 | 4.77 | | 40 | 0.0709 | 3.92 |
| | 60 | 0.0309 | 0.69 | | 60 | 0.0191 | 3.23 |
| $\overline{\frac{\partial u}{\partial y}}^x D/U$ | 20 | 0.1564 | - | $\overline{\frac{\partial v}{\partial y}}^x D/U$ | 20 | 0.2260 | - |
| | 40 | 0.0255 | 2.62 | | 40 | 0.0305 | 2.89 |
| | 60 | 0.0102 | 2.26 | | 60 | 0.0158 | 1.62 |
| $\overline{(u)_A}^{xy}/U$ | 20 | 0.0277 | - | $\overline{(v)_A}^{xy}/U$ | 20 | 0.0352 | - |
| | 40 | 0.0033 | 3.06 | | 40 | 0.0022 | 4.00 |
| | 60 | 0.0014 | 2.17 | | 60 | 0.00095 | 2.07 |
| $\overline{(u)_B}^{xy}/U$ | 20 | 0.0123 | - | $\overline{(v)_B}^{xy}/U$ | 20 | 0.0142 | - |
| | 40 | 0.0027 | 2.21 | | 40 | 0.0025 | 2.48 |
| | 60 | 0.0013 | 1.81 | | 60 | 0.0010 | 2.28 |
| $\overline{(u)_{1-4}}^{xy}/U$ | 20 | 0.0208 | - | $\overline{(v)_{1-4}}^{xy}/U$ | 20 | 0.0088 | - |
| | 40 | 0.0013 | 4.02 | | 40 | 0.00061 | 3.84 |
| | 60 | 0.00058 | 1.94 | | 60 | 0.00028 | 1.97 |

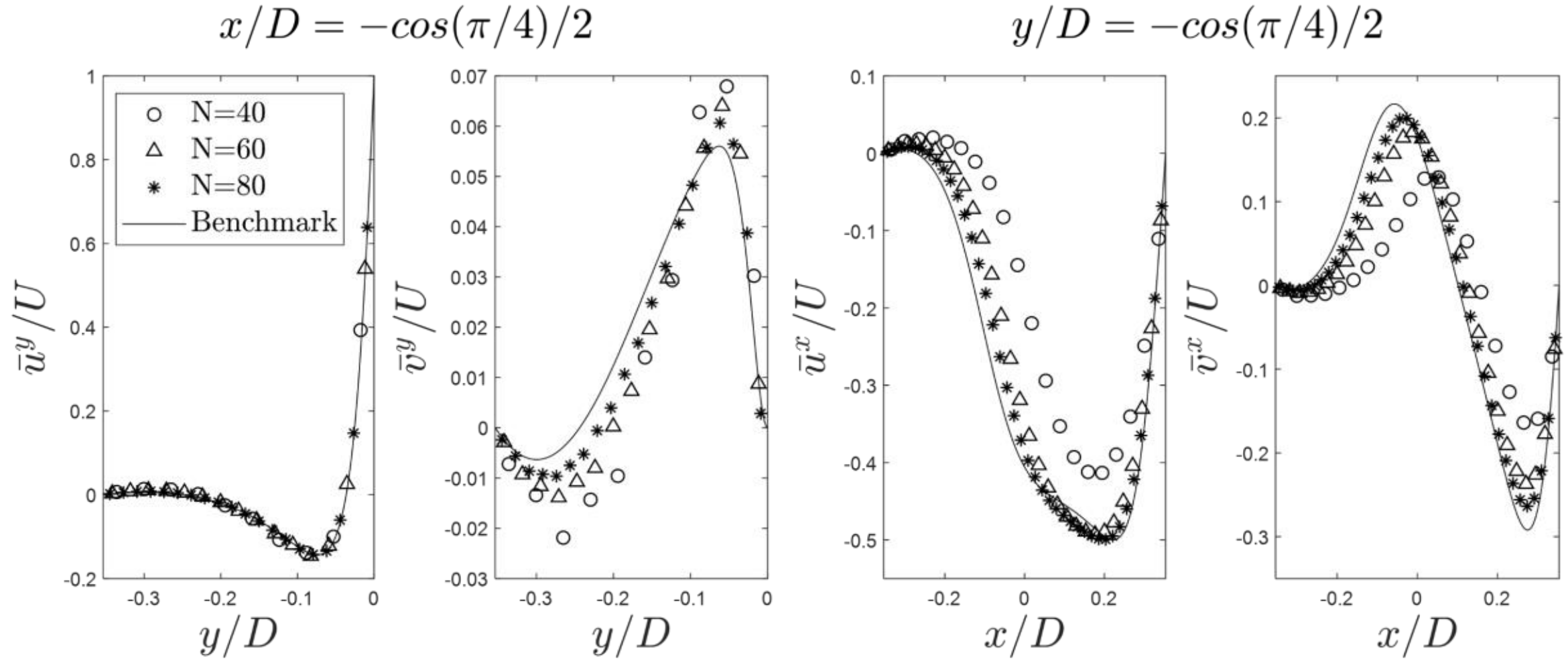

*Figure 16*

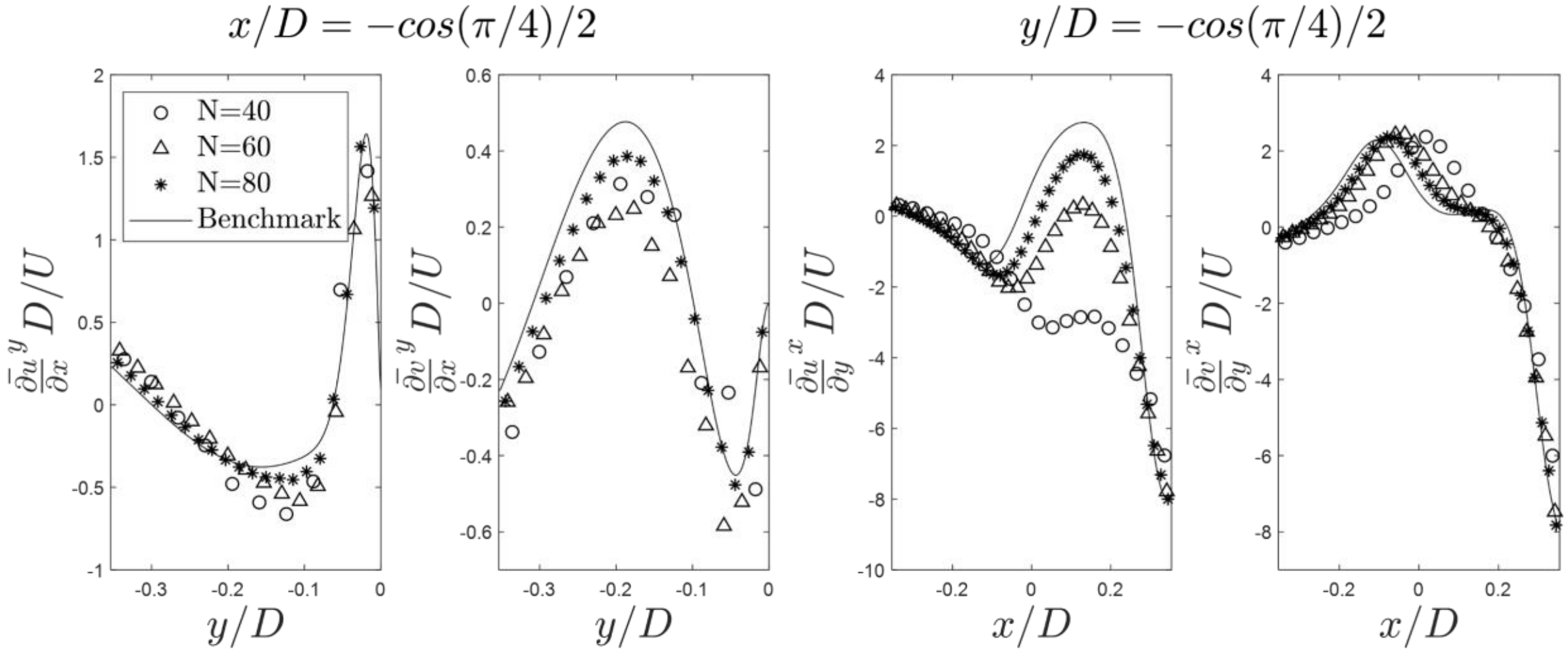

*Figure 17*

Cells A

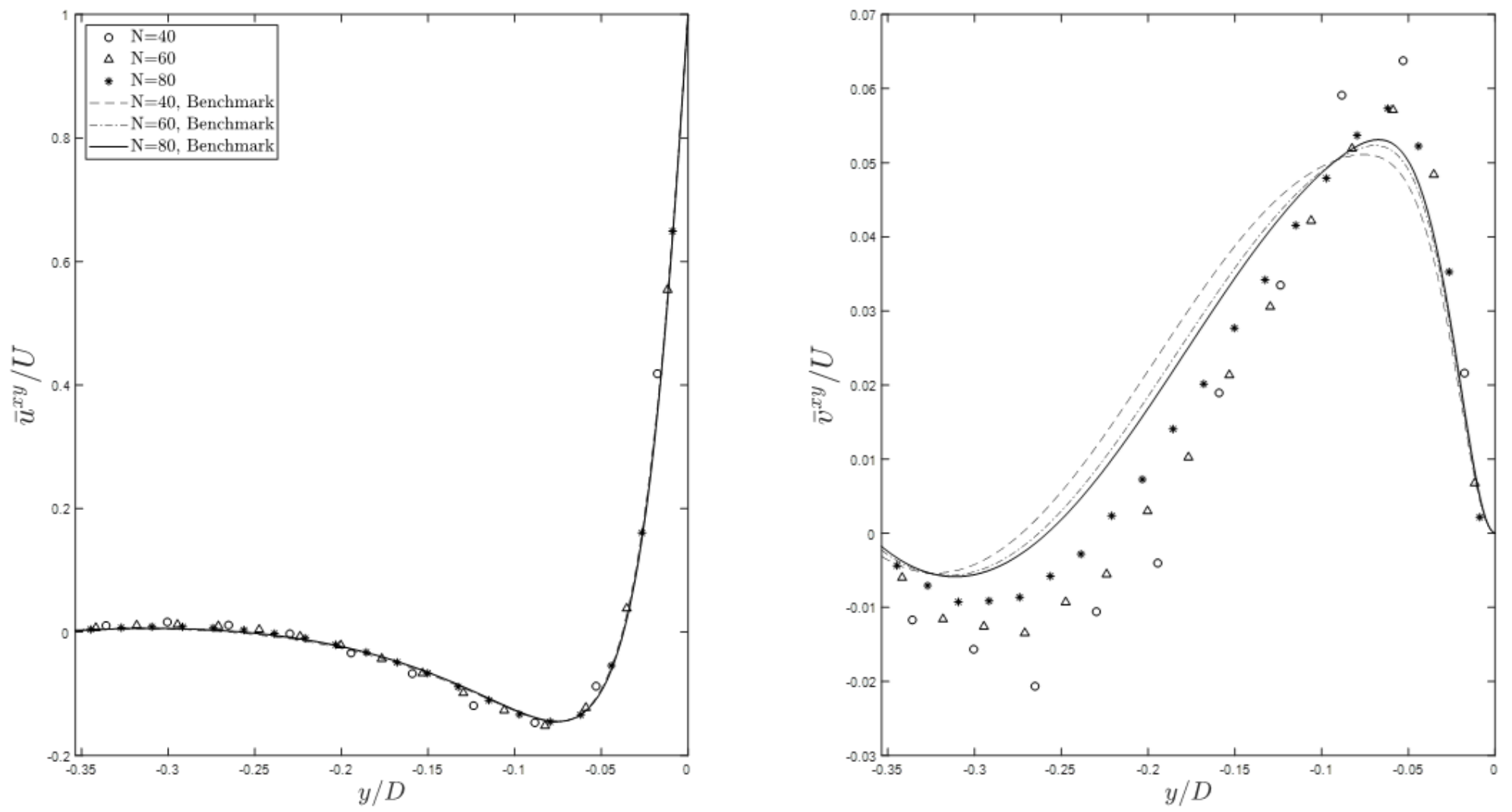

*Figure 18*

Cells B

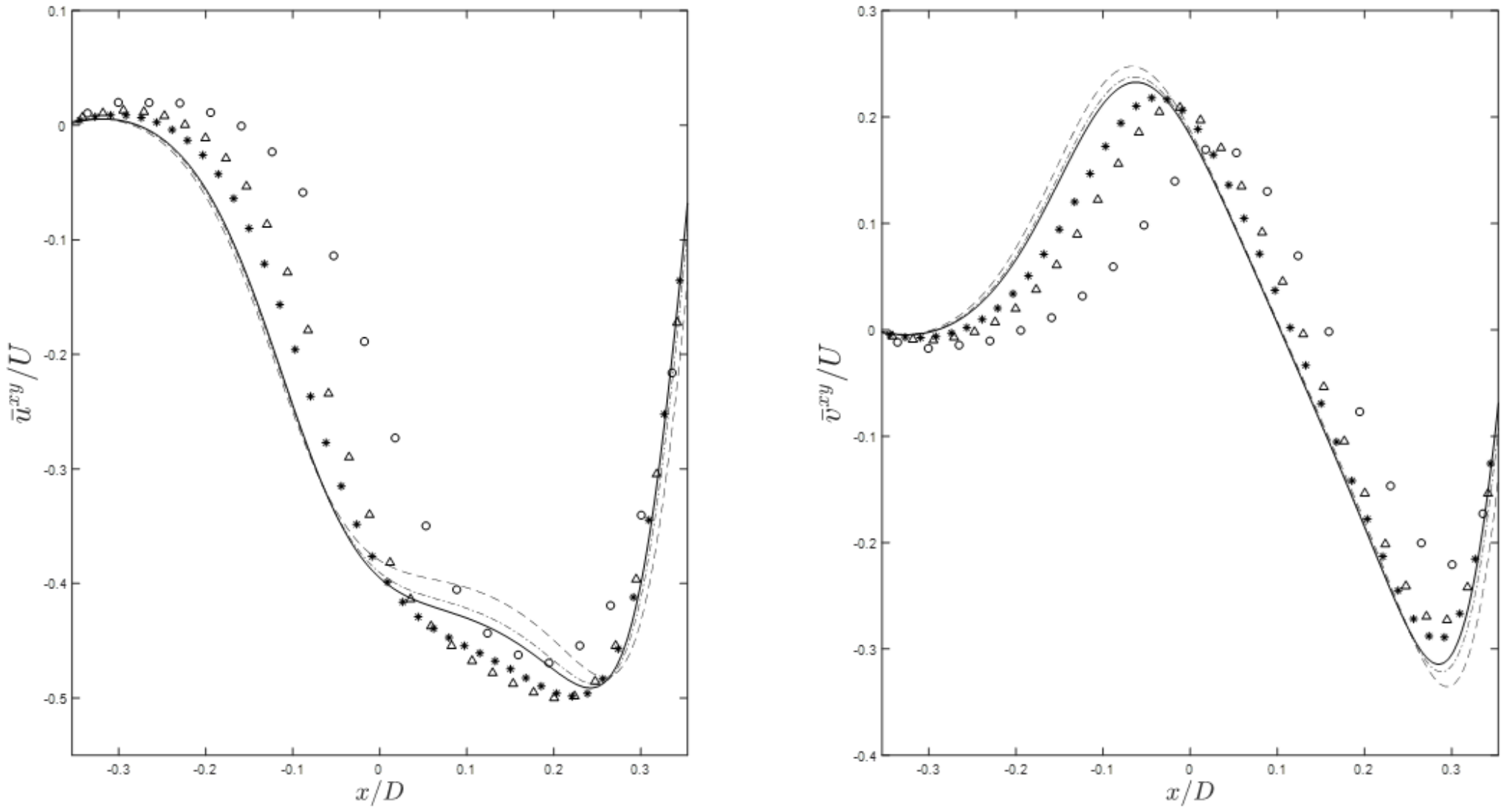

*Figure 19*

*Table 3*

| Velocity | Mesh ($N$) | Cell 1 | Cell 2 | Cell 3 | Cell 4 |
|---|---|---|---|---|---|
| $\bar{u}^{xy}/U$ | 40 | 0.0058 | 0.0029 | -0.0865 | -0.0744 |
| | Benchmark | 0.0022 | 0.0022 | -0.0916 | -0.0862 |
| | 60 | 0.0020 | 0.0028 | -0.0634 | -0.0584 |
| | Benchmark | 0.0016 | 0.0016 | -0.0615 | -0.0596 |
| | 80 | 0.0016 | 0.0015 | -0.0478 | -0.0460 |
| | Benchmark | 0.0013 | 0.0012 | -0.0470 | -0.0464 |
| $\bar{v}^{xy}/U$ | 40 | -0.0061 | -0.0026 | -0.0786 | -0.0511 |
| | Benchmark | -0.0024 | -0.0020 | -0.0977 | -0.0750 |
| | 60 | -0.0020 | -0.0026 | -0.0605 | -0.0475 |
| | Benchmark | -0.0017 | -0.0015 | -0.0641 | -0.0543 |
| | 80 | -0.0018 | -0.0013 | -0.0463 | -0.0400 |
| | Benchmark | -0.0013 | -0.0011 | -0.0484 | -0.0432 |

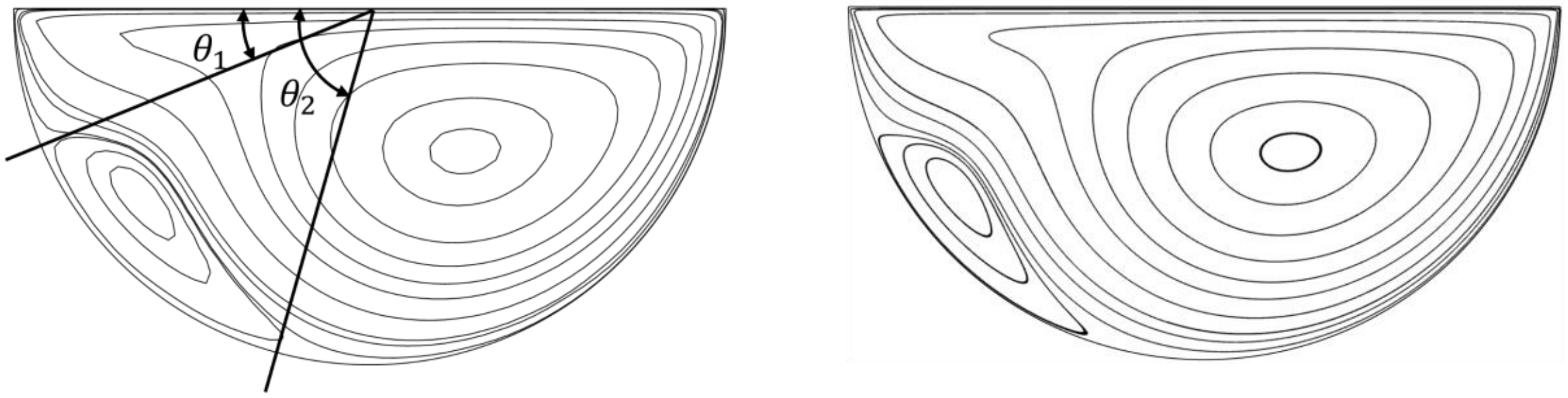

*Figure 20*

*Table 4*

| Variable | Mesh ($N$) | Error | Order of accuracy | Variable | Mesh ($N$) | Error | Order of accuracy |
|---|---|---|---|---|---|---|---|
| $\bar{u}^y/U$ | 40 | 0.0115 | - | $\bar{v}^y/U$ | 40 | 0.0131 | - |
| | 60 | 0.0054 | 1.85 | | 60 | 0.0073 | 1.44 |
| | 80 | 0.0025 | 2.64 | | 80 | 0.0043 | 1.84 |
| $\bar{u}^x/U$ | 40 | 0.1088 | - | $\bar{v}^x/U$ | 40 | 0.0800 | - |
| | 60 | 0.0402 | 2.46 | | 60 | 0.0392 | 1.76 |
| | 80 | 0.0210 | 2.26 | | 80 | 0.0206 | 2.23 |
| $\overline{\frac{\partial u}{\partial x}}^y D/U$ | 40 | 0.1860 | - | $\overline{\frac{\partial v}{\partial x}}^y D/U$ | 40 | 0.1581 | - |
| | 60 | 0.1250 | 0.98 | | 60 | 0.1730 | -0.22 |
| | 80 | 0.0540 | 2.92 | | 80 | 0.0599 | 3.69 |
| $\overline{\frac{\partial u}{\partial y}}^x D/U$ | 40 | 2.3331 | - | $\overline{\frac{\partial v}{\partial y}}^x D/U$ | 40 | 0.9141 | - |
| | 60 | 1.2816 | 1.48 | | 60 | 0.5119 | 1.43 |
| | 80 | 0.6382 | 2.42 | | 80 | 0.2887 | 1.99 |
| $\overline{(u)_A}^{xy}/U$ | 40 | 0.0148 | - | $\overline{(v)_A}^{xy}/U$ | 40 | 0.0149 | - |
| | 60 | 0.0069 | 1.89 | | 60 | 0.0089 | 1.27 |
| | 80 | 0.0031 | 2.74 | | 80 | 0.0048 | 2.12 |
| $\overline{(u)_B}^{xy}/U$ | 40 | 0.0788 | - | $\overline{(v)_B}^{xy}/U$ | 40 | 0.0867 | - |
| | 60 | 0.0398 | 1.68 | | 60 | 0.0405 | 1.88 |
| | 80 | 0.0211 | 2.21 | | 80 | 0.0211 | 2.27 |
| $\overline{(u)_{1-4}}^{xy}/U$ | 40 | 0.0053 | - | $\overline{(v)_{1-4}}^{xy}/U$ | 40 | 0.0118 | - |
| | 60 | 0.0012 | 3.70 | | 60 | 0.0029 | 3.43 |
| | 80 | 0.00048 | 3.16 | | 80 | 0.0015 | 2.35 |

Table 5

| | Current work | Benchmark | Relative error |
|---|---|---|---|
| $\psi_{min}/DU$ | -0.0766 | -0.0779 | 1.7% |
| $(x/D, y/D)$ | (0.1321, -0.1945) | (0.1214, -0.2030) | (8.8%, 4.2%) |
| $\theta_1$ | 22.74° | 21.42° | 6.2% |
| $\theta_2$ | 74.74° | 71.49° | 4.5% |

# 4. Summary and Discussion

Compared with other finite-volume based complex boundary treatment methods, our method has the following advantages:

1. We use cut Cartesian grids. Such grids can be generated with minimal programming efforts and stored with minimal computer memories.
2. We do not approximate complex boundaries with piecewise linear segments. Such approximations are practiced by almost all other methods except those using coordinate transformation and may introduce uncertain errors into the numerical calculation.
3. Instead of using polynomial interpolation or reconstruction, we employ only calculus and finite differences to obtain second-order accurate fluxes and gradients on irregular faces. Polynomial interpolation/reconstruction usually entails solving a large number of equations to determine the coefficients in reconstruction polynomials. And these equations themselves involve very significant number of complicated integrals that have to be evaluated, although only once typically for a given mesh. Our method also requires evaluation of some integrals once, but they are much simpler. During each iteration or time step, instead of solving a system of equations for coefficients of reconstruction polynomials, our method only needs calculation of certain finite differences, which demands fewer computing resources, especially considering the fact that such solving systems of equations and calculation of finite differences must be carried out for every mesh cell.
4. Some methods have to distinguish quite many different cell configurations and apply different treatments accordingly. Our method only identifies four different situations: rectangular cells, irregular boundary cells, solitary cells, and twin cells.
5. Our method allows enormous flexibility when it is implemented. Indeed, the most essential component of this method is the PKP theorem, and almost all other technical details can be customized. In this sense the current method is more a general framework than a particular numerical scheme.

6. The finite difference equations resulted from this method are of familiar sparse diagonal systems of equations, which invite many well-established solving methods.

Some disadvantages of the current method include:

1. Our way of generating the cut cell mesh is quite rigid, less flexible than many other methods.
2. Because of this specific way of mesh generation, some grid cells have very large aspect ratios.
3. There seems no easy way to extend the current method to three-dimensional flows except for the extrusion of two-dimensional geometry along the third-dimension case.
4. It might be difficult to realize fourth or higher order schemes in the current framework. This point is discussed in more details now.

Extension of the current method to fourth order accuracy involves using fifth-order accurate PKP theorem to calculate sliding averages of derivatives of flow field variables along cell faces (refer to Equations (56) and (62)). One then has to explicitly evaluate the fourth-order terms on the right side of the theorem (see Equation (24)) which includes many new derivatives as well as using at least third-order accurate finite differences to approximate the two derivatives in the second term on the right side of Equation (24). Although theoretically possible, such a method is complicated. It seems, however, relatively straightforward to extend the current method to third-order accuracy as such an extension still uses the original fourth order PKP theorem. Of course, all face fluxes and gradients then must be evaluated with third-order accurate equations. Such a scheme will be very useful in for example LES of turbulent flows in complex geometries.